\documentclass[graybox]{svmult}

\usepackage{type1cm,bm}        
\usepackage{comment}
\usepackage{makeidx}         
\usepackage{graphicx}        
\usepackage{multicol}        
\usepackage[bottom]{footmisc}

\usepackage{mathrsfs}

\usepackage{newtxtext}       %
\usepackage[varvw]{newtxmath}       

\usepackage{hyperref}
\hypersetup{
    colorlinks,
    citecolor=blue,
    filecolor=blue,
    linkcolor=blue,
    urlcolor=blue
}

\spnewtheorem{myexample}{Example}{\bfseries}{\itshape} 
\newcommand{\bx}{{\bm x}}
\newcommand{\bn}{{\bm n}}
\newcommand{\bv}{{\bm v}}

\newcommand{\Uad}{{\mathscr U_\mathsf{ad}}}
\newcommand{\Mad}{{\mathscr M_\mathsf{ad}}}
\newcommand{\Xad}{{X_\mathsf{ad}}}
\newcommand{\R}{\mathbb{R}}
\newcommand{\N}{\mathbb{N}}

\newcommand{\tilv}{\tilde{v}}
\DeclareMathOperator*{\conv}{conv}
\DeclareMathOperator*{\argmin}{arg\,min}

\makeindex             

\begin{document}

\title*{Multiobjective Optimization of Non-Smooth PDE-Constrained Problems}
\author{M. Bernreuther, M. Dellnitz, B. Gebken, G. M\"uller, S. Peitz, K. Sonntag and S. Volkwein}
\institute{Michael Dellnitz, Bennet Gebken \& Konstantin Sonntag \at Department of Mathematics, Paderborn University, Germany, \email{dellnitz@upb.de / bgebken@math.upb.de / konstantin.sonntag@upb.de}
\and Georg M\"uller \at University of Heidelberg, D-69120 Heidelberg, Germany, \email{georg.mueller@uni-heidelberg.de},
\and Sebastian Peitz \at Department of Computer Science, Paderborn University, Germany, \email{sebastian.peitz@uni-paderborn.de}
\and Marco Bernreuther \& Stefan Volkwein \at University of Konstanz, D-78457 Konstanz, Germany, \email{stefan.volkwein@uni-konstanz.de}}
%
%
\maketitle

\abstract{
Multiobjective optimization plays an increasingly important role in modern applications, where several criteria are often of equal importance. The task in multiobjective optimization and multiobjective optimal control is therefore to compute the set of optimal compromises (the Pareto set) between the conflicting objectives. The advances in algorithms and the increasing interest in Pareto-optimal solutions have led to a wide range of new applications related to optimal and feedback control -- potentially with non-smoothness both on the level of the objectives or in the system dynamics. This results in new challenges such as dealing with expensive models (e.g., governed by partial differential equations (PDEs)) and developing dedicated algorithms handling the non-smoothness. Since in contrast to single-objective optimization, the Pareto set generally consists of an infinite number of solutions, the computational effort can quickly become challenging, which is particularly problematic when the objectives are costly to evaluate or when a solution has to be presented very quickly.
This article gives an overview of recent developments in the field of multiobjective optimization of non-smooth PDE-constrained problems. In particular we report on the advances achieved within Project 2 ``Multiobjective Optimization of Non-Smooth PDE-Constrained Problems -- Switches, State Constraints and Model Order Reduction'' of the DFG Priority Programm 1962 ``Non-smooth and Complementarity-based Distributed Parameter Systems: Simulation and Hierarchical Optimization''.
}

\section{Introduction}
\label{Section:1}
Multicriteria decision making problems are ubiquitous in all areas of daily life. For instance, we want to produce goods as cheaply as possible while maintaining a high quality. In autonomous driving, we want to reach a destination as quickly as possible while consuming a minimal amount of energy. In the same way, political decisions have to take economical, societal and ecological criteria into account simultaneously.
These examples give rise to the Pareto principle, named after the italian scientist Vilfredo Pareto (1848 -- 1923). Since the objectives are generally conflicting, it is not possible to satisfy them all at the same time. Instead, we are forced to accept compromises. According to the Pareto principle, a compromise can be called optimal only if it is impossible to further improve all criteria simultaneously. The solution to such a \emph{multicriteria decision making problem} (also \emph{multiobjective optimization problem (MOP)} is the set of optimal compromises, also referred to as the \emph{Pareto set}. It's image under the objective functions is then called the \emph{Pareto front}.
In order to support decision makers, the task in multiobjective optimization is to provide approximations of the entire set. This way, we gain a lot of insight into possible trade-off solutions and can reach decisions in a much more informed way.

Due to the fact that the solution is a set -- in contrast to isolated minimizers for single-criteria problems -- the numerical computation of solutions to MOPs is in general significantly more challenging. This fact becomes even more significant if the underlying objective functions are expensive to evaluate numerically, or if the problem possesses less regularity than usual (e.g., we can only assume Lip\-schitz continuity).
We here want an overview of recent developments in the field of non-smooth multiobjective optimization of problems whose objective functions require the evaluation of partial differential equations (PDEs). The results covered here also include the advances achieved within Project 2 ``Multiobjective Optimization of Non-Smooth PDE-Constrained Problems -- Switches, State Constraints and Model Order Reduction'' of the DFG Priority Programm 1962 ``Non-smooth and Complementarity-based Distributed Parameter Systems: Simulation and Hierarchical Optimization''.
In particular -- after introducing the basics of multiobjective optimization (Section \ref{Section:2}) -- we are going to address the following sub-topics:
\begin{itemize}
    \item Efficient algorithms for the solution of non-smooth MOPs, both without and with PDE constraints (Section \ref{Section:3}).
    \item How to incorporate surrogate models for PDEs and treat the resulting inexactness (Section \ref{Section:4}).
    \item Some applications of (non-smooth) multiobjective optimization in machine learning (Section \ref{Section:5}).
    \item How to exploit non-smoothness (in the form of mixed-integer control) in order to significantly reduce the surrogate modeling effort of smooth PDE-constrained control problems (Section \ref{Section:6}).
\end{itemize}

\section{Basic concepts}
\label{Section:2}

\subsection{Multiobjective Optimization}
\label{Section:2.1}

In the following we introduce the basics of multiobjective optimization (MO). For further details we refer to \cite{Ehr05,Mie12}. MO is concerned with the simultaneous optimization of multiple objective functions. Let $\Xad$ be a convex, closed set and let $f_i:\Xad\to\R$, $i\in\{1,\dots,k\}$, be the different objectives, which may be grouped into the objective vector $f:\Xad\rightarrow \R^k$. We first have to clarify what it means for a point to be ``optimal'' in this case. For $k=1$, i.e., for scalar optimization, a point $x\in\Xad$ is \emph{optimal} if there is no point $y\in\Xad$ with $f(y)<f(x)$. For $k>1$, this concept can be generalized as in the following definition, where $<$ and $\leq$ for elements of $\R^k$ are meant component-wise.

\begin{definition} \label{def:pareto_optimality}
    A point $x\in\Xad$ is \emph{weakly Pareto optimal}, if there is no $y\in\Xad$ with $f(y)<f(x)$. It is \emph{Pareto optimal}, if there is no $y \in X$ with $f(y) \leq f(x)$ and $f(y) \neq f(x)$. The set of Pareto optimal points is the \emph{Pareto set}, denoted by $P$. The image $f(P) \subseteq \R^k$ of the Pareto set is the \emph{Pareto front}.
\end{definition}

The problem of finding the Pareto set of $f$ is called a \emph{multiobjective optimization problem} (MOP) and is denoted as
\begin{align}
    \label{eq:general_MOP}
    \min_{x\in\Xad} f(x) \quad \text{or} \quad \min_{x\in\Xad} \ 
    \begin{pmatrix}
        f_1(x) \\
        \vdots \\
        f_k(x)
    \end{pmatrix}.
\end{align}
While the definition of weak Pareto optimality is formally identical to classical notion of optimality, there is a crucial difference that is caused by the properties of ``$<$'' and ``$\leq$'' on $\R^k$ for $k > 1$: For $k = 1$ and $a, b \in \R^k$, we have either $a < b$, $a = b$ or $a > b$. But for $k > 1$, e.g., $k = 2$, we may have $a = (1,0)^\top$ and $b = (0,1)^\top$, such that $a \nless b$, $a \neq b$ and $a \ngtr b$. Thus, in this case, neither of the points $a$ and $b$ is better than the other one. In light of Definition \ref{def:pareto_optimality}, this means that the Pareto set $P$ typically contains more than one point. Furthermore, it implies that the optimal value, i.e., the Pareto front, is not a unique value as in the scalar case, but a whole set of optimal values. The following example visualizes this behavior.

\begin{myexample}
    \label{example:paraboloids}
    For $\Xad= \R^n$, $k = 2$ and $c^1, c^2 \in \R^n$, consider 
    \begin{align*}
        f : \R^2 \rightarrow \R^2, \quad x \mapsto 
        \begin{pmatrix}
            {\| x - c^1 \|}_2^2 \\[0.5mm]
            {\| x - c^2 \|}_2^2
        \end{pmatrix},
    \end{align*}
    where $\|\cdot\|_2$ stands for the Euclidean norm. Then a point $x \in \R^n$ is Pareto optimal if there is no other point $y$ that has at least the same distances to $c^1$ and $c^2$ as $x$, while being closer to one of the points. It is easy to see that this is equivalent to $x$ lying in the convex hull $\conv(\{c^1,c^2\})$ of $c^1$ and $c^2$. Figure~{\em\ref{fig:example_paraboloids}} shows a visualization of this example for $n = 2$, $c^1 = (0,0)^\top$ and $c^2 = (1,1/2)^\top$.
    \begin{figure}
        \centering
        \parbox[b]{0.49\textwidth}{
            \centering 
            \includegraphics[width=0.4\textwidth]{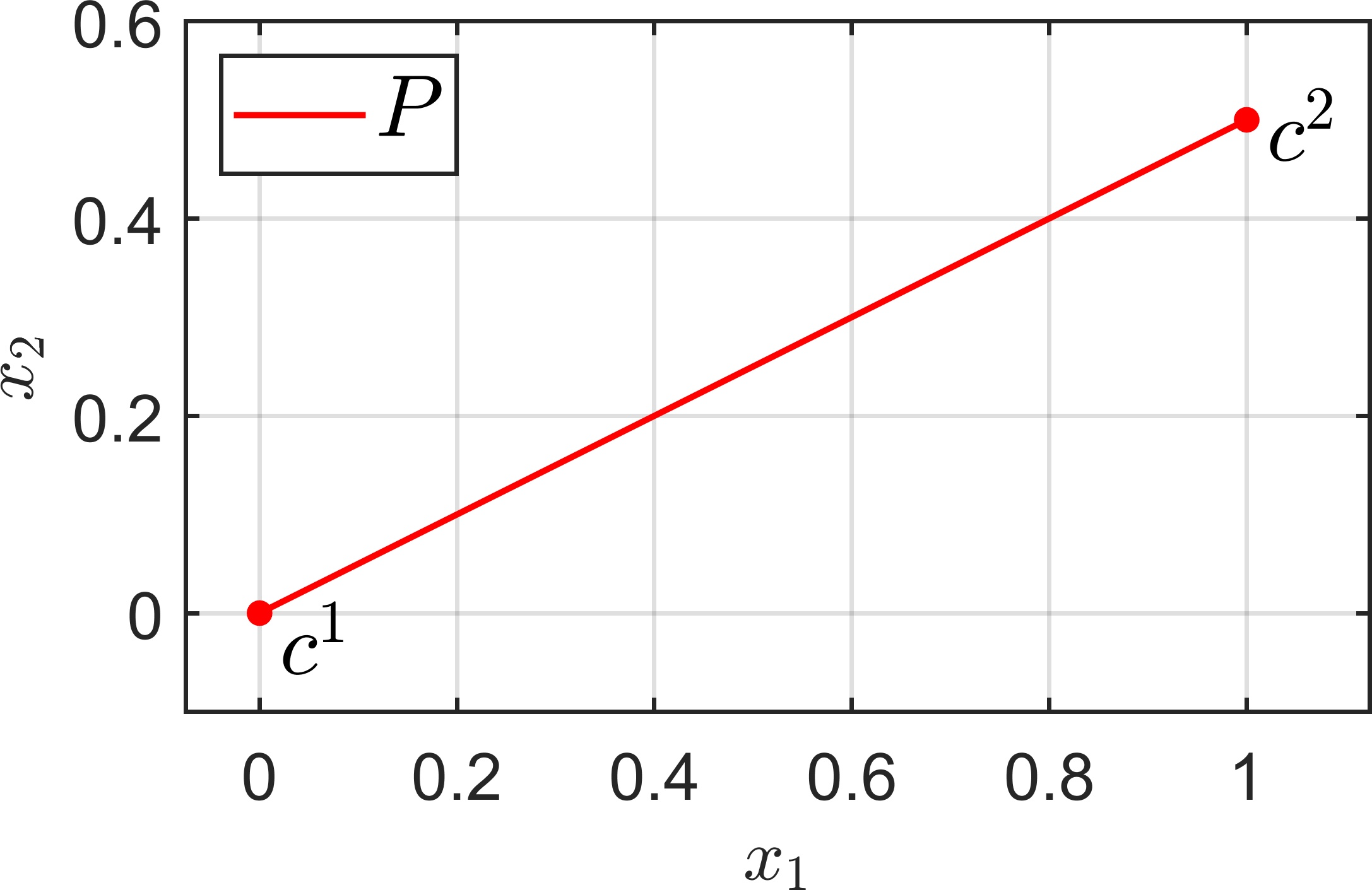}\\
            (a)
		}
        \parbox[b]{0.49\textwidth}{
            \centering 
            \includegraphics[width=0.4\textwidth]{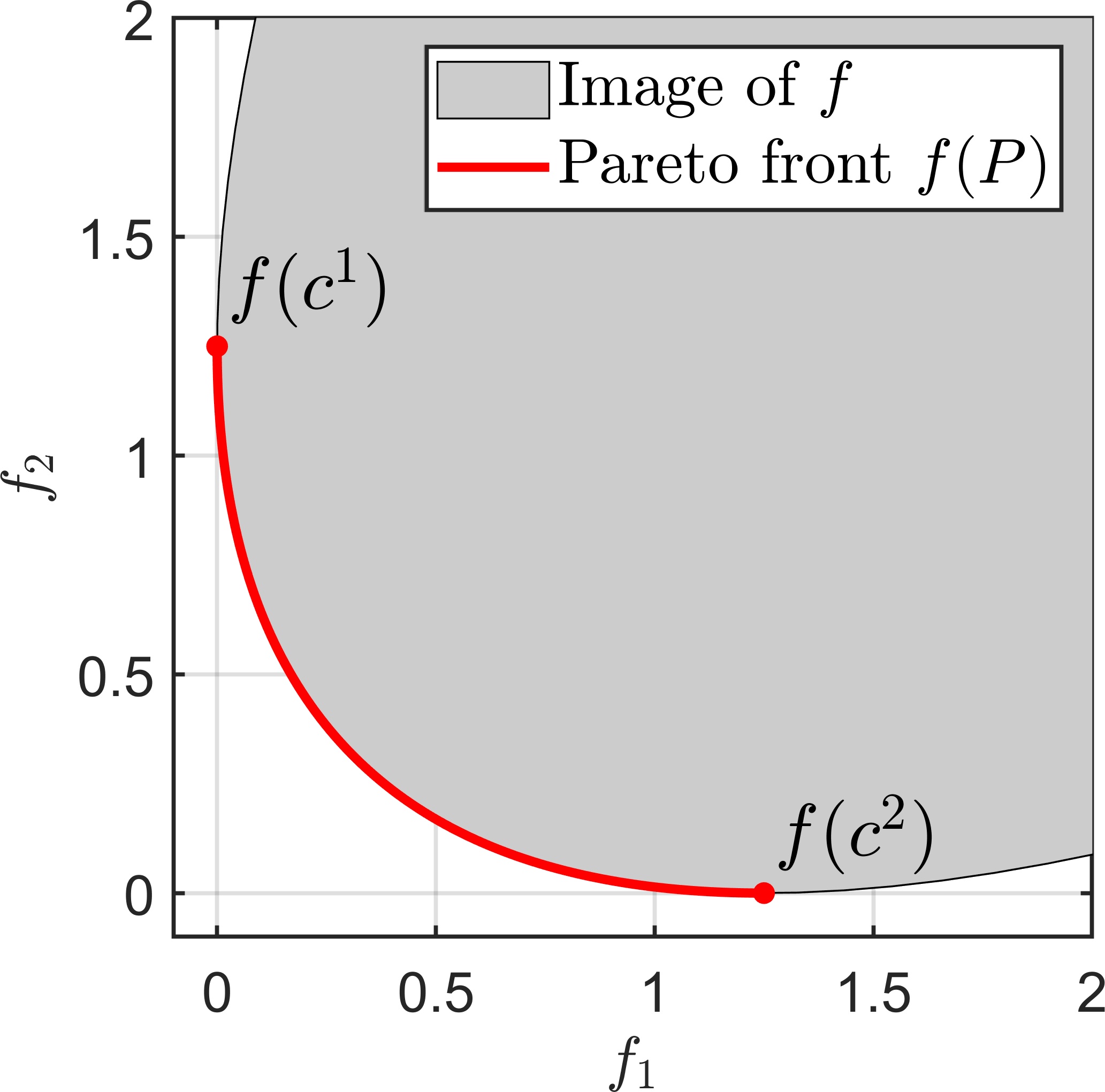}\\
            (b)
		}
        \caption{Example \ref{example:paraboloids} for $n=2$, $c^1 = (0,0)^\top$ and $c^2 = (1,1/2)^\top$: (a) The Pareto set $P$. (b) The Pareto front $f(P)$ and the image $f(\Xad)$ of $f$.}
        \label{fig:example_paraboloids}
    \end{figure}
\end{myexample}

As in the scalar case, numerical solution of MOPs requires optimality conditions. In case $\Xad= \R^n$ and $f$ is differentiable, a first-order necessary condition for weak Pareto optimality is the existence of a vanishing convex combination of the gradients of the objectives; cf., e.g., \cite{Ehr05,Mie12}. More formally, if $\bar x$ is weakly Pareto optimal, then
\begin{align}
    \label{eq:smooth_opt_condition}
    0 \in \conv(\{ \nabla f_1(\bar x), \dots, \nabla f_k(\bar x) \}),
\end{align}
i.e., 
\begin{align}
    \label{eq:smooth_KKT}
    \exists\bar\alpha \in (\R^{\geq 0})^k:\sum_{i=1}^k\bar\alpha_i = 1\text{ and }\sum_{i=1}^k\bar\alpha_i\nabla f_i(\bar x) = 0.
\end{align}
Analogously to critical points in the scalar case, points that satisfy \eqref{eq:smooth_opt_condition} are called \emph{Pareto critical}.


\subsection{PDE-Constrained Multiobjective Optimization}
\label{Section:2.2}

When the objective functions $f_i$ and their gradients $\nabla f_i$ are expensive to evaluate -- e.g., since the constraints are given by PDEs -- the computational time for the numerical solution of MOPs can quickly become very or even too large. In that case surrogate models offer a promising tool to reduce the computational effort significantly. Examples are dimensional reduction techniques are the Reduced Basis (RB) method (cf., e.g., \cite{Haa17,Hesthaven,QMN16}) or Proper Orthogonal Decomposition (POD) (cf., e.g., \cite{Vol17}).

\begin{myexample}
    \label{Example:PDEs-1}
    Let us present two examples. We suppose that $\Omega\subset\R^d$, $d\in\{1,2,3\}$, is a bounded domain with Lipschitz-continuous boundary $\Gamma=\partial\Omega$. The outward normal vector is denoted by $\bn$. Moreover, let $H$ and $V$ denote the standard separable Hilbert spaces $L^2(\Omega)$ and $H^1(\Omega)$, respectively, endowed by the usual inner products
    \begin{align*}
        {\langle\varphi,\psi \rangle}_H=\int_\Omega\varphi\psi\,\mathrm d\bx,\quad{\langle\varphi,\psi\rangle}_V=\int_\Omega\varphi\psi+\nabla\varphi\cdot\nabla\psi\,\mathrm d\bx
    \end{align*}
    and their associated induced norms. We write $V'$ for the dual space of $V$. For more details on Lebesgue and Sobolev spaces we refer to {\em\cite{Eva10}}, for instance.
    \begin{enumerate}
        \item [\em 1)] For $T>0$ we set $Q=(0,T)\times\Omega$, $\Sigma=(0,T)\times\Gamma$. Recall the Hilbert space
        \begin{align*}
             W(0,T)=\big\{\varphi\in L^2(0,T;V)\,\big|\,\varphi_t\in L^2(0,T;V')\big\}
        \end{align*}
        endowed with the common inner product
        \begin{align*}
            {\langle\varphi,\phi\rangle}_{W(0,T)}=\int_0^T{\langle \varphi(t),\phi(t)\rangle}_V+{\langle\varphi_t(t),\phi_t(t)\rangle}_{V'}\,\mathrm dt    
        \end{align*}
        and the induced norm; compare {\em\cite[pp.~472-479]{DL00}}, for instance. Here, the expression $\varphi(t)$ stands for the function $\varphi(t,\cdot)$ considered as a function in $\Omega$ only.\hfill\\
        It is well-known that $W(0,T)$ is continuously embedded into $C([0,T];H)$, the space of continuous functions from $[0, T]$ to $H$. Let $(\mathscr H,\langle\cdot\,,\cdot\rangle_\mathscr H)$ a given (observation) Hilbert space endowed with the induced norm $\|\cdot\|_\mathscr H$. We utilize the (control) Hilbert space $\mathscr U=L^2(0,T;\R^m)$ supplied with the usual topology and set $\mathcal Y=W(0,T)$. Obviously, we have $L^2(0,T;V)\hookrightarrow \mathscr Y$. Now we consider the following (infinite-dimensional) MOP:
        \begin{subequations}
            \label{MOPparabolic}
	        \begin{equation}
                \label{minipro-1}
                \min J(y,u)=\left(
		        \begin{aligned}
                    J_1(y,u)\\
                    J_2(y,u)
		        \end{aligned}
                \right)=\frac{1}{2}\,\left(
                \begin{array}{c}
                    {\|\mathcal Cy-y^\mathsf d\|}_\mathscr H^2\\[1mm] 
                    \displaystyle\sum_{i=1}^m\int_0^T{\|u_i(t)\|}_2^2\,\mathrm dt
                \end{array}
                \right) 
	        \end{equation}
	        subject to $(y,u)\in \mathscr Y\times\Uad$ weakly solves the parabolic boundary value problem
	        \begin{equation}
                \label{parpde}		
			    \begin{aligned}
			         y_t(t,\bx)-\Delta y(t,\bx)+\bv(\bx)\cdot\nabla y(t,\bx)&=g(t,\bx)&&\text{for }(t,\bx)\in Q,\\
                    \frac{\partial y}{\partial \bn}(t,\bx)&= \sum_{i=1}^mu_i(t)\xi_i(\bx)&&\text{for }(t,\bx)\in\Sigma,\\
                    y(0,\bx)&=y_\circ(\bx)&&\text{for }\bx\in\Omega.   
			    \end{aligned}
	        \end{equation}
        \end{subequations}
        Here, $\mathcal C:L^2(0,T;V)\to\mathscr H$ is a linear, bounded (observation) operator, $\Uad$ a closed convex set, $g\in L^2(0,T;H)$ an inhomogeinity, $y_\circ\in H$ an initial condition, and $\xi_1,\ldots,\xi_m\in L^2(\Gamma)$ are control shape functions. Furthermore, the convection field $\bv$ belongs to $H^1(\Omega;\R^d)$ and satisfies $\nabla\cdot\bv=0$. It follows from standard arguments {\em\cite{DL00,Eva10}} that \eqref{parpde} has a unique weak solution $y\in\mathscr Y$ for every $u\in\Uad$. Thus, we define the (affin-linear) control-to-state operator $\mathcal S:\Uad\to W(0,T)$, where $y=\mathcal S(u)$ weakly solves \eqref{parpde} for given $u\in\Uad$. This allows us to define the reduced cost functional $f(u)=J(\mathcal S(u),u)\in\R^2$ for every $u\in\Uad$ and to consider -- instead of \eqref{MOPparabolic} the reduced MOP (cf. \eqref{eq:general_MOP})
        \begin{equation}
            \label{PDEEx1-RedMOP}
            \min f(u)=\left(
		        \begin{array}{c}
                    f_1(u)\\
                    \vdots\\
                    f_k(u)
		        \end{array}
                \right)\quad\text{s.t.}\quad u\in\Uad.
        \end{equation}
        We refer the reader to {\em\cite{Ban17,Ban21,BBV16,BBV17,BFGV20,Mak20,Spu19}}, where \eqref{PDEEx1-RedMOP} is studied analytically and numerically.
        \item [\em 2)] Let us consider the following non-convex MOP
        \begin{subequations}
            \label{MOPelliptic}
	        \begin{equation}
                \label{minipro-2}
                \min J(y,\mu)=\left(
		        \begin{array}{c}
                    J_1(y,\mu)\\
                    \vdots\\
                    J_k(y,\mu)
		        \end{array}
                \right)=\frac{1}{2}\,\left(
                \begin{array}{c}
                    \int_\Omega |y-y_1^\mathsf d|^2\,\mathrm d\bx\\ 
		            \vdots\\
		            \int_\Omega |y - y_{k-1}^\mathsf d|^2 \, \mathrm d\bx \\[1mm] 
		            \textstyle\sum\limits_{j=1}^m|\mu_j - \mu^\mathsf d_j|^2
                \end{array}
                \right) 
	        \end{equation}
	        subject to $(y,\mu)\in V\times\Mad$ weakly solves the semi-linear elliptic boundary value problem
	        \begin{equation}
                \label{semipde}		
			    \begin{aligned}
			        -\Delta y(\bx)+y(\bx)+y(\bx)^3&= g(\bx)+ \sum_{i = 1}^m\mu_i \xi_i(\bx)&&\text{for }\bx\in\Omega,\\
                    \frac{\partial y}{\partial \bn}(\bx) + y(\bx)&= g_\Gamma(\bx)&&\text{for }\bx\in\Gamma,   
			    \end{aligned}
	        \end{equation}
        \end{subequations}
        where $y \in V$ is the state variable and $\mu\in\Mad= [\mu_\mathsf a,\mu_\mathsf b]$ the parameter. We suppose that $g_\Gamma \in L^r(\Gamma)$ with $r>d-1$, $\xi_1,...,\xi_m,g\in H$, $\mu^\mathsf d=(\mu_j^\mathsf d)\in\R^m$ and $y_1^\mathsf d,...,y_{k-1}^\mathsf d \in H$.\hfill\\
        It is well-known that \eqref{semipde} has a unique weak solution $y=y(\mu)\in V$ for every $\mu\in\mathscr M_\mathsf{ad}$. Thus, the non-linear solution operator $\mathcal S:\Mad\to V$, where $y=\mathcal S(\mu)$ is the unique weak solution to \eqref{semipde} is well-defined. We define the reduced cost functional $f(\mu)=J(\mathcal S(\mu),\mu)\in\R^k$  for every $\mu\in\mathscr M_\mathsf{ad}$ and express \eqref{MOPelliptic} as
        \begin{equation}
            \label{PDEEx2-RedMOP}
            \min f(\mu)=\left(
		        \begin{array}{c}
                    f_1(\mu)\\
                    \vdots\\
                    f_k(\mu)
		        \end{array}
                \right)\quad\text{s.t.}\quad\mu\in\Mad
        \end{equation}
        which is of the form \eqref{eq:general_MOP}. In particular, we have that $\Mad$ is even compact. Problem \eqref{PDEEx2-RedMOP} is studied analytically and numerically in {\em\cite{BGDPV22,BGRV21,BMV22a,IUV17,Rei20}}, for instance.
    \end{enumerate}
\end{myexample}

Notice, that evaluating the reduced cost functional $f$ at a given point $x$ ($x=u$ for \eqref{PDEEx1-RedMOP} and $x=\mu$ for \eqref{PDEEx2-RedMOP}) requires the computation of the solution operator $\mathcal S$ at $x$, i.e., the solution of a PDE. This is usually computationally expensive so that we are interested in techniques to accelerate the evaluation process. Here, reduced-order methods are very suitable to speed-up the computational time while still being accurate enough. The latter is ensured by a-posteriori error estimates. 

As explained in Section~\ref{Section:2.1} weak Pareto optimal points can be characterized by first-order necessary optimality conditions. Let us assume that $X$ is a Hilbert space and $\Xad\subset X$ a convex, closed subset. Taking the convexity constraints into account, it follows that if $\bar x$ is weakly Pareto optimal then there exist $\bar \alpha \in (\R^{\geq 0})^k$ with $\sum_{i = 1}^k \bar\alpha_i = 1$ and
\begin{align}
    \label{OptCondPDE}
    \sum_{i = 1}^k\bar\alpha_i\,{\langle\nabla f_i(\bar x),x-\bar x\rangle}_X\ge0\quad\text{for all }x\in\Xad,
\end{align}
where $\nabla f_i(\bar x)\in X$ denotes the gradient of the $i$-th reduced coct functional at $\bar x$. Utilizing the adjoint approach \cite{HPUU08} the gradients $\nabla f_i$ can be characterized by adjoint variables $p_i$ for $i=1,\ldots,k$. However, the computation of the adjoint (or dual) variables requires $k$ additional PDE solves so that we also apply reduced-order methods to approximate the adjoint variables.

\begin{myexample}
    \label{Example:PDEs-2}
    Let us come back to Example~{\em\ref{Example:PDEs-1}} and present the formulas for the gradients $\nabla f_i(x)$.
    \begin{enumerate}
        \item [\em 1)] Suppose that an arbitrarily $u\in\Uad$ the state $y=\mathcal S(u)$ is given as the weak solution to \eqref{parpde}. Let $p\in\mathscr Y$ be the weak solution to the backward equation
	    \begin{equation}
            \label{parpdedual}		
			\begin{aligned}
			    \quad\hspace{-7mm}-p_t(t,\bx)-\Delta p(t,\bx)-(\bv\cdot\nabla p)(t,\bx)&=\big(\mathcal C^*(y^\mathsf d-\mathcal Cy)\big)(t,\bx)&&\hspace{-2mm}\text{for }(t,\bx)\in Q,\\
                \frac{\partial p}{\partial \bn}(t,\bx)+\big(\bv(\bx)\cdot\bn\big)p(t,\bx)&=0&&\hspace{-2mm}\text{for }(t,\bx)\in\Sigma,\\
                p(T,\bx)&=0&&\hspace{-2mm}\text{for }\bx\in\Omega,  
			\end{aligned}
	    \end{equation}
        where $\mathcal C^*:\mathscr H\to L^2(0,T;V)$ denotes the (Hilbert space) adjoint of $\mathcal C$ satisfying
        \begin{align*}
             {\langle\mathcal C^*\phi,\varphi\rangle}_{L^2(0,T;V)}={\langle\phi,\mathcal C\varphi\rangle}_\mathscr H\quad\text{for all }\phi\in\mathscr H,\varphi\in L^2(0,T;V).
        \end{align*}
        Now we have
        \begin{align*}
            \nabla f_1(u)=\left(
            \begin{array}{c}
                 -\int_\Gamma\xi_1(\bx) p(\cdot\,,\bx)\,\mathrm d\bx\\
                 \vdots\\
                 -\int_\Gamma\xi_m(\bx) p(\cdot\,,\bx)\,\mathrm d\bx
            \end{array}
            \right)\in\mathscr U,\quad\nabla f_2(u)=\left(
            \begin{array}{c}
                 u_1\\
                 \vdots\\
                 u_m
            \end{array}
            \right)\in\mathscr U.
        \end{align*}
        Consequently, to get the gradient $\nabla f_1(u)$ we have to solve the state equation \eqref{parpde} and the adjoint equation \eqref{parpdedual}, i.e., we need two PDE solves.
        \item [\em 2)] Let $\mu\in\Mad$ be given arbitrarily and $y=\mathcal S(\mu)$ solve \eqref{semipde} weakly. By $p_i\in W(0,T)$, $i=1,\ldots,k-1$, we denote the adjoint variables solving
        \begin{align}
            \label{semipdeDual}		
			\begin{aligned}
			    -\Delta p_i(\bx)+p_i(\bx)+3y(\bx)^2p_i(\bx)&=y_i^\mathsf d(\bx)-y(\bx)&&\text{for }\bx\in\Omega,\\
                \frac{\partial p_i}{\partial \bn}(\bx) + p_i(\bx)&=0&&\text{for }\bx\in\Gamma. 
			\end{aligned}
	    \end{align}
        Now the gradients are given as ($i=1,\ldots,k-1$)
        \begin{align*}
            \nabla f_i(\mu)=\left(
            \begin{array}{c}
                 -\int_\Omega\xi_1(\bx) p_i(\cdot\,,\bx)\,\mathrm d\bx\\
                 \vdots\\
                 -\int_\Omega\xi_m(\bx) p_i(\cdot\,,\bx)\,\mathrm d\bx
            \end{array}
            \right)\in\R^m,\quad\nabla f_k(\mu)=\left(
            \begin{array}{c}
                 \mu_1-\mu_1^\mathsf d\\
                 \vdots\\
                 \mu_m-\mu_m^\mathsf d
            \end{array}
            \right)\in\R^m.
        \end{align*}
        so that we need to solve the state equation \eqref{semipde} and $k-1$ adjoint equations \eqref{semipdeDual} to compute the gradients $\nabla f_i$, $i=1,\ldots,k-1$. In this case we have to compute (weak) solutions to $k$ PDEs.
    \end{enumerate}
\end{myexample}

\subsection{Non-smoothness}
\label{Section:2.3}

In our setting there are two ways in which non-smoothness may occur: The objective vector $f$ of the MOP may be a non-smooth function and/or the PDE may contain non-smooth terms. We will first discuss non-smoothness in the objective (for $X = \R^n$). For general results on non-smooth multiobjective optimization we refer the reader to \cite{Mak02,MEK14,MKW16}.

\paragraph{Non-smoothness in the objective vector}

For the objective vector $f$, we follow the standard assumption in non-smooth optimization that its components $f_i$, $i \in \{1,\dots,k\}$, are only locally Lipschitz continuous. Since the gradients of the objectives $f_i$ generally do not exist in this case, the optimality condition \eqref{eq:smooth_opt_condition} (and all methods based on it) cannot be used. Instead, concepts from non-smooth analysis \cite{C1990} may be used. More precisely, the gradient of $f_i$ may be replaced by the \emph{Clarke subdifferential}
\begin{equation} \label{eq:def_clarke_subdiff}
    \begin{aligned}
        \partial f_i(x) := \conv &\left( \left\{ \xi \in \R^n : \exists (x^j)_j \in \R^n \setminus \Omega \text{ with } \lim_{j \rightarrow \infty} x^j = x \text{ and } \right. \right. \\
        &\left. \left. \lim_{j \rightarrow \infty} \nabla f_i(x^j) = \xi \right\} \right),
    \end{aligned}
\end{equation}
where $\Omega_i \subseteq \R^n$ is the set of points in which $f_i$ is not differentiable. In the scalar case, a necessary condition for optimality is that $0 \in \partial f(\bar x)$ for an optimal solution $\bar x$. For the multiobjective case, this condition may be generalized as follows \cite{MEK14}: If $\bar x$ is weakly Pareto optimal and does not lie on the boundary of $\Xad$, then
\begin{align}
    \label{eq:nonsmooth_opt_condition}
    0 \in \conv \left( \left\{ \partial f_1(\bar x) \cup \dots \cup \partial f_k(\bar x) \right\} \right).
\end{align}

While the Clarke subdifferential can be used analogously to the smooth gradient in many theoretical results, it is much more difficult to use in practice: Since the Clarke subdifferential $\partial f_i(x)$ only captures the non-smoothness of $f_i$ if $x$ is a non-smooth point of $f_i$, and since $\Omega_i$ is a null set (by Rademacher's theorem), we cannot use it to treat non-smoothness in practice. Instead, the so-called \emph{(Goldstein)} $\varepsilon$\emph{-subdifferential} is used (see \cite{G1977})
\begin{align}
    \partial_\varepsilon f_i(x) := \conv \left( \bigcup_{y \in B_\varepsilon(x)} \partial f_i(y) \right),
\end{align}
where $\varepsilon \geq 0$ and $B_\varepsilon(x) := \{ y \in \R^n : \| x - y \|_2 \leq \varepsilon \}$, may be used. For $\varepsilon > 0$, the $\varepsilon$-subdifferential may be interpreted as a stabilized Clarke subdifferential. In Section \ref{Section:3.1} we will show how practical solution methods for non-smooth MOPs can be derived from it. 

\begin{remark}
    Let us mention two related recent works on the treatment of non-smooth objectives in the context of PDE constraints.
    \begin{enumerate}
        \item [1)] For a convex, closet subset $\Mad\in\R^m$ the following minimization problem is considered in \cite{BDPV17}:
        \begin{align*}
	       \min J(y,\mu)=\frac{1}{2}\left( 
            \begin{array}{c}
                {\|y-y_1^\mathsf d\|}_{L^2(\Omega)}^2\\[1mm]
                {\|y-y_2^\mathsf d\|}_{L^1(\Omega)}\\[1mm]
                \sum\limits_{i=1}^m|\mu_i|^2
            \end{array}
            \right)
        \end{align*}
        subject to $(y,\mu)\in V\times\Mad$ weakly solves the semilinear PDE constraints
        \begin{align*}
	       -\Delta y(\bx)+y(\bx)^3=\sum_{i = 1}^m\mu_i \xi_i(\bx)\text{ for }\bx\in\Omega,\quad\frac{\partial y}{\partial \bn}(\bx)=0\text{ for }\bx\in\Gamma,
        \end{align*}
        Here, $\Omega$, $\Gamma$, $H$, $V$ are as in Example~\ref{Example:PDEs-1}-2), $y_1^\mathsf d,y_2^\mathsf d\in H$.
        \item [2)] We refer to the work \cite{AB23} that provides a comprehensive study of the non-monotone forward-backward splitting method for solving
        \begin{align*}
             \min_{x\in X} f(x):=f_1(x)+f_2(x),
        \end{align*}
        where $f_1$ is a continuously Fr\'echet-differentiable function (possibly non-convex), and $f_2$ is a convex function whose proximal operator
        \begin{align*}
            \mathrm{prox}_{f_2}(x):=\argmin_{\tilde x\in X}\left(f_2(\tilde x)+\frac{1}{2}\,{\|\tilde x-x\|}_X^2\right)
        \end{align*}
        is assumed to be explicitly computable.
    \end{enumerate}    
\end{remark}

\paragraph{Non-smoothness in the PDE}

The non-smoothness can occur in several ways. On one hand there can occur non-smooth non-linearities. We have considered here the $\max$ function both in elliptic and parabolic PDEs (see, e.g., \cite{BGRV21,BMG22b,Ber19,BV22,Ber19,BMV22b,CM21,Rei20}, respectively). On the other hand, mixed-integer optimal control problems also leads to non-smoothness in the PDE; cf. \cite{Jae19,JV20}.

\begin{myexample}
    \label{Ex:MaxPDE}
    Here we just briefly recall the elliptic PDE with the max function. Let $\Omega\subset\R^d$, $d\in\{1,2,3\}$, be a bounded domain with Lipschitz-continuous boundary $\Gamma=\partial\Omega$. We set $H=L^2(\Omega)$ and $V=H^1_0(\Omega)$. Suppose that $\Mad=[\mu_\mathsf a,\mu_\mathsf b]\subset\R^m$ with $\mu_\mathsf a\le\mu_\mathsf b$ component-wise. Then, for given control $u\in\mathscr U=H$ and parameter $\mu\in\Mad$we consider the parametrized equation
\begin{align}
    \label{MaxPDE1}
    \begin{aligned}
        -c(\mu)\Delta y(\bx)+a(\mu)\max\big\{0,y(\bx)\big\}&=u(\bx)&&\text{for }\bx\in\Omega,\\ y(\bx)&=0&&\text{for }\bx\in\Gamma,    
    \end{aligned}
\end{align}
where $c:\Mad\to\R^{>0}$ and $a:\Mad\to\R^{\ge0}$ are Lipschitz continuous. A function $y$ is called a weak solution to \eqref{MaxPDE1} provided $y\in V$ and for all $\varphi\in V$ we have
\begin{align}
    \label{MaxPDE2}
    \int_\Omega c(\mu)\nabla y(\bx)\cdot\nabla \varphi(\bx)+a(\mu)\max\big\{0,y(\bx)\big\}\varphi(\bx)\,\mathrm d\bx=\int_\Omega u(\bx)\varphi(\bx)\,\mathrm d\bx.
\end{align}
To solve \eqref{MaxPDE2} numerically we apply a finite element (FE) Galerkin discretization. Let $\{\varphi_i\}_{i=1}^N\subset V$ be $N$ linear independent FE functions and $V^N=\mathrm{span}\,\{\varphi_1,\ldots,\varphi_N\}\subset V$. Then, we look for the FE function
\begin{align}
    \label{MaxPDE2a}
     y_\mathsf{fe}(\bx)=\sum_{j=1}^N\mathrm y_j\varphi_j(\bx)\quad\text{for }\bx\in\Omega
\end{align} 
with unknowns $\mathrm y_1,\ldots,\mathrm y_N\in\R$ satisfying
\begin{align}
    \label{MaxPDE3}
    \begin{aligned}
        &\int_\Omega c(\mu)\nabla y_\mathsf{fe}(\bx)\cdot\nabla \varphi_i(\bx)+a(\mu)\max\big\{0,y_\mathsf{fe}(\bx)\big\}\varphi_i(\bx)\,\mathrm d\bx\\
        &\hspace{20mm}=\int_\Omega u(\bx)\varphi_i(\bx)\,\mathrm d\bx    
    \end{aligned}
\end{align}
for $i=1,\ldots,N$. Inserting \eqref{MaxPDE2a} into \eqref{MaxPDE3} problem \eqref{MaxPDE3} can be expressed as an $N$-dimensional nonlinear algebraic system that is solved -- due to the non-smoothness -- efficiently by semi-smooth Newton method in {\em\cite{Ber19,BMV22c,BMV22b}}.
\end{myexample}


\section{Algorithms for non-smooth multiobjective optimization}
\label{Section:3}

In this section, we will present solution methods for non-smooth MOPs. The first one is a descent method for general non-smooth MOPs \cite{GP21a,Geb22a} and the second one is a method that is tailored to a specific class of non-smooth MOPs that arises in the context of regularization \cite{BGP22,GBP23}. 

\subsection{Descent method for non-smooth MOPs}
\label{Section:3.1}

Here, we will consider general MOPs \eqref{eq:general_MOP} for the case where the $f_i$, $i \in \{1,\dots,k\}$, are only locally Lipschitz continuous. We will fist consider the case where $X = \R^n$ and then briefly discuss the generalization to the case where $X$ is an arbitrary Hilbert space.

The idea of descent methods is to generate a sequence $x^j \in \R^n$, starting in some $x^1 \in \R^n$, such that $f(x^{j+1}) < f(x^j)$ for all $j \in \N$. Given some $x \in \R^n$, we achieve this by first computing a direction $v \in \R^n$ in which all objectives decrease and then computing a step length $t > 0$ such that $f(x + t v) < f(x)$. In theory, convex analysis shows that such a direction may be computed as
\begin{align}
    \label{eq:def_bar_v}
    \bar{v} := \argmin\big\{{\| \xi \|}_2^2\,\big|\,\xi \in -\partial_\varepsilon^\cup f(x)\big\}
\end{align}
where $\partial_\varepsilon^\cup f(x) := \conv(\bigcup_{i = 1}^k \partial_\varepsilon f_i(x))$. The direction $\bar{v}$ is, in a sense, the ``steepest'' descent direction and satisfies 
\begin{align}
    \label{eq:bar_v_armijo}
    f_i(x + t \bar{v}) \leq f_i(x) - t \| \bar{v} \|_2^2 \quad \text{for all }t \in \left( 0, \frac{\varepsilon}{\| \bar{v} \|_2} \right],\,i \in \{1,\dots,k\},
\end{align}
which can be used to compute a step length via a standard backtracking line search. However, since the entire $\varepsilon$-subdifferential required for the computation of $\bar{v}$ is generally not available in practice, a method based on solving \eqref{eq:def_bar_v} has little use. 

Instead, we compute an approximation $\bar{v}$ based on an approximation of $\partial_\varepsilon^\cup f(x)$ that only requires us to be able to evaluate an arbitrary subgradient at every $x \in \R^n$. The idea is to start with an approximation $W = \{ \xi_1, \dots, \xi_k \} \subseteq \partial_\varepsilon^\cup f(x)$ that only contains one subgradient per objective $f_i$ and to then iteratively add new subgradients to $W$ until an acceptable descent direction is found.  To this end, for $W \subseteq \partial_\varepsilon^\cup f(x)$ let
\begin{align*}
    \tilv := \argmin \big\{{\| \xi \|}_2^2\,\big|\,\xi \in -\conv(W)\big\}.
\end{align*}
Based on \eqref{eq:bar_v_armijo}, we say that $\tilv$ is an \emph{acceptable} descent direction if 
\begin{align*}
    f_i \left( x + \frac{\varepsilon}{\| \tilv \|_2} \bar{v} \right) \leq f_i(x) - c \varepsilon\,{\| \bar{v} \|}_2 \quad \forall i \in \{1,\dots,k\}
\end{align*}
for some fixed $c \in (0,1)$. To find a new subgradient that improves the approximation of $\partial_\varepsilon^\cup f(x)$ in case $\tilv$ is not acceptable for some $i \in \{1,\dots,k\}$, the mean value theorem from non-smooth analysis implies that there are $t' \in (0,\varepsilon / \| \tilv \|_2)$ and $\xi' \in \partial f_i(x + t' \tilv)$ with
\begin{align*}
    \langle \xi', \tilv \rangle > -c\,{\| \tilv \|}_2^2.
\end{align*}
Convex analysis shows that such a $\xi'$ cannot be contained in $\conv(W)$, such that adding $\xi'$ to $W$ improves the approximation of $\partial_\varepsilon^\cup f(x)$. In \cite{GP21a} it was shown that iteratively adding subgradients to $W$ in this way leads to an acceptable direction $\tilv$ after finitely many iterations. Using a backtracking line search yields a simple descent method for non-smooth MOPs. In \cite{Geb22a} it was shown that when dynamically reducing $\varepsilon$ to zero, then the method converges to Pareto critical points of $f$.
\begin{figure}
    \centering
    \parbox[b]{0.49\textwidth}{
        \centering 
        \includegraphics[width=0.45\textwidth]{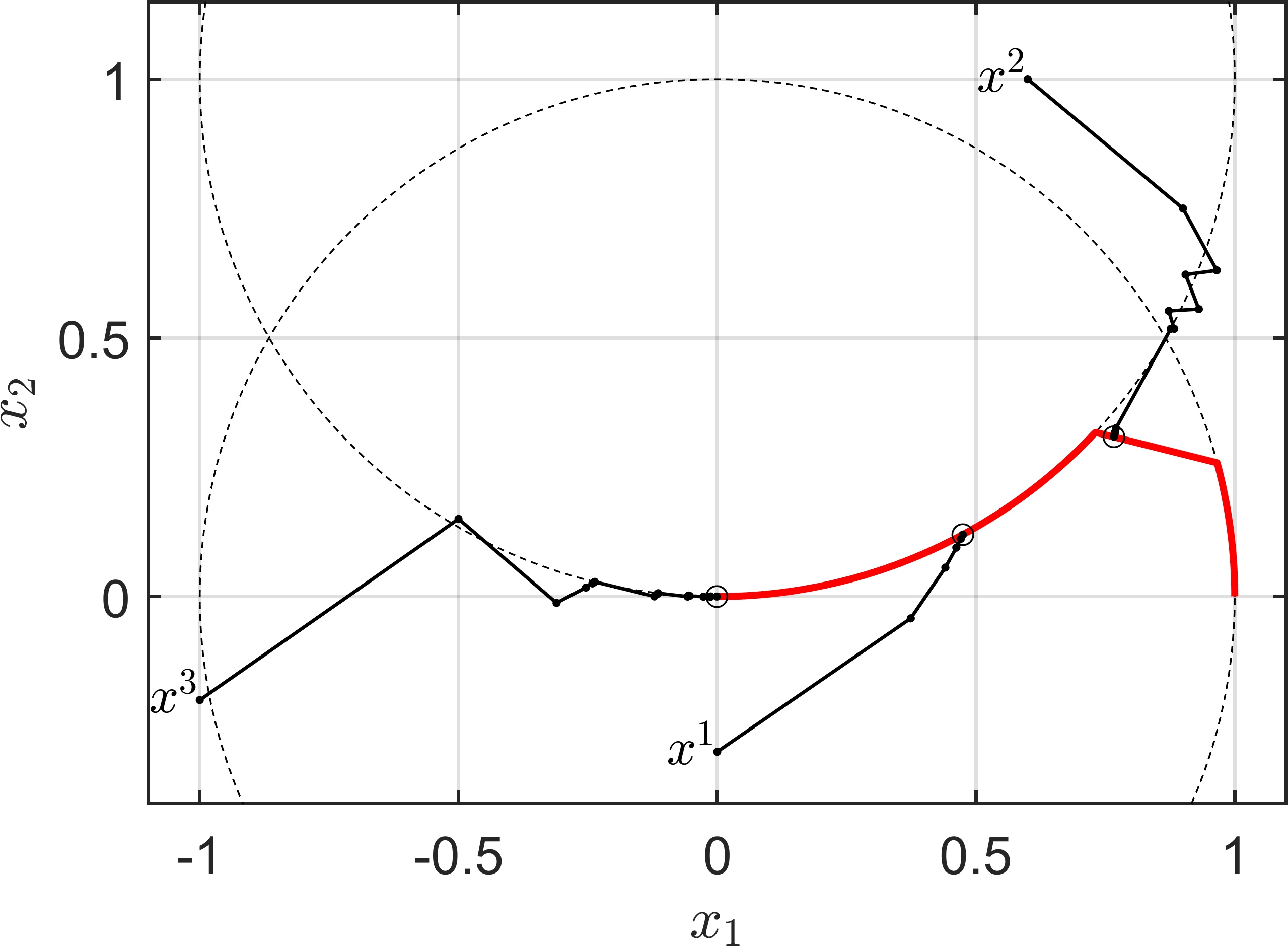}\\
        (a)
	}
    \parbox[b]{0.49\textwidth}{
        \centering 
        \includegraphics[width=0.45\textwidth]{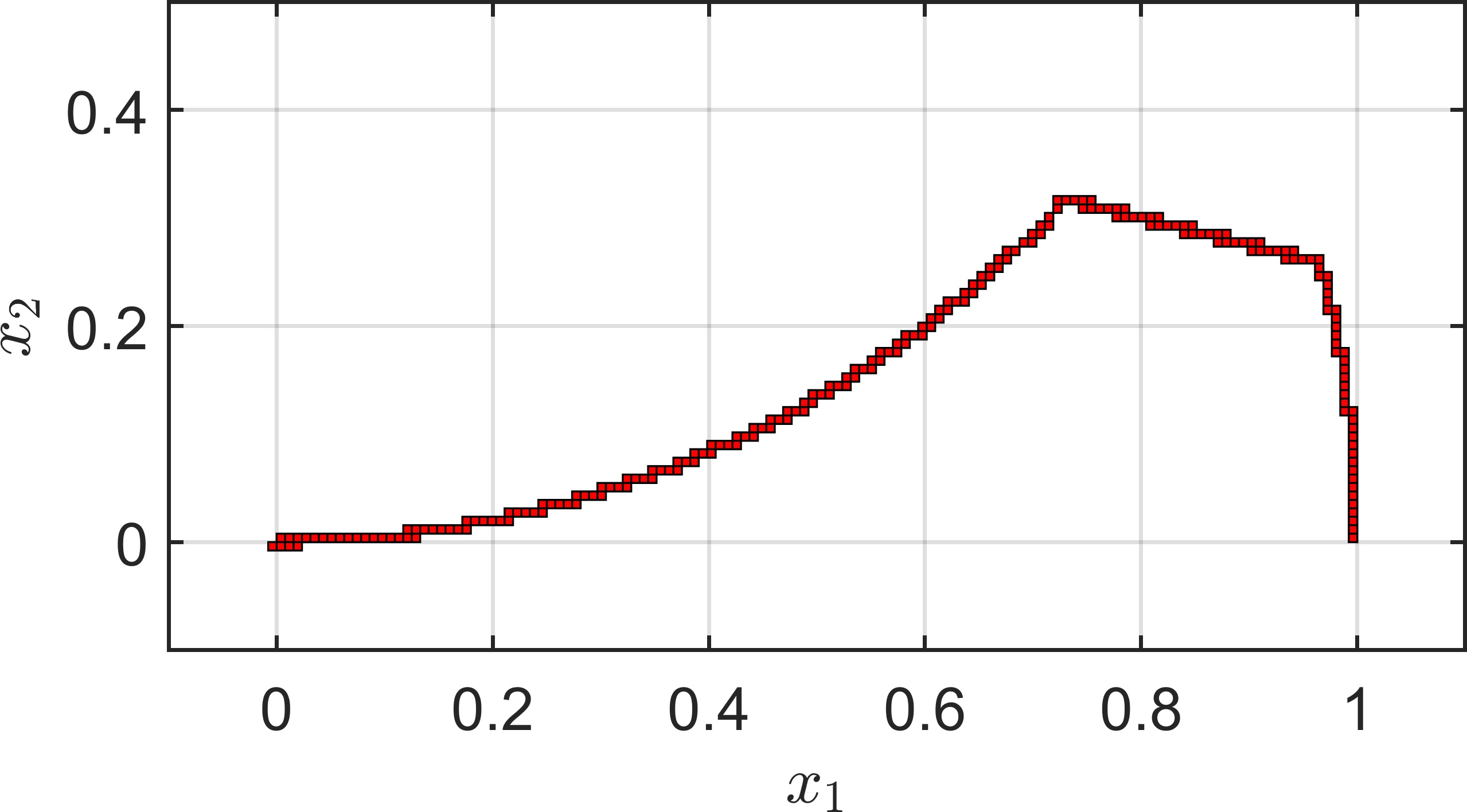}\\
        (b)
	}
    \caption{(a) The descent method from Section \ref{Section:3.1} applied to different starting points for the problem M30 from \cite{MKW2014}. (b) A box covering of the Pareto set. (Taken from \cite{GP21a} License: CC BY 4.0).}
    \label{fig:example_descent_method}
\end{figure}  
Figure \ref{fig:example_descent_method}(a) shows the behavior of the method in an example. Numerical experiments showed that in terms of subgradient evaluations, the method is competitive to the \emph{multiobjective proximal bundle method} \cite{MKW16}. Finally, to compute the entire Pareto set, the \emph{subdivision techniques} from \cite{DSH2005} can be used. This yields a method that is able to compute box coverings of the Pareto set, as shown in Figure \ref{fig:example_descent_method}(b).


The descent method can be generalized to MOPs \eqref{eq:general_MOP} with objective functions defined on a general Hilbert space $X = \mathcal{H}$ which is infinite dimensional. When going from the Euclidean space $\R^n$ to $\mathcal{H}$, we have to take care since some objects lose properties which get used in order to show that the descent method computes Pareto critical solutions. This is in particular the case for the compactness of the multiobjective $\varepsilon$-Goldstein subdifferential. Therefore, we cannot directly use the $\varepsilon$-Goldstein subdifferential but instead use its (weak$^*$) closure
\begin{align*}
    F_{\varepsilon}(x) = \overline{\conv} \left( \bigcup_{i = 1}^k \overline{\partial_{\varepsilon} f_i(y)} \right) \subset \mathcal{H}^*.
\end{align*}
If the functions $f_i$ are Lipschitz continuous on $B_{\overline{\varepsilon}}(x)$, the set $F_{\varepsilon}(x)$ is nonempty, convex and weakly compact for all $\varepsilon \in \left[0, \overline{\varepsilon}\right)$. Keeping this in mind, the proof of convergence for the descent method in infinite dimensions works almost analogously to the finite dimensional setting.

The generalization to Hilbert spaces allows the treatment of non-smooth optimal control problems. We consider an obstacle problem on a two-dimensional domain $\Omega = (-1,1)^2$ with two objective functions. Given an obstacle $\psi \in H^1(\Omega)$, define the set of admissible displacements
\begin{align*}
    K \coloneqq \left\lbrace y \in H_0^1(\Omega) : y \le \psi \text{ a.e. on } \Omega \right\rbrace.
\end{align*}
Set $V = H_0^1(\Omega)$ and $H = L^2(\Omega)$. The constraining obstacle problem can be described in variational form as
\begin{align*}
    \text{Find } y \in K : {\langle Ay - b, v - y \rangle}_{V', V} \ge 0 \text{ for all } v \in K,
\end{align*}
where $A:V \to V'$ is a linear, bounded and coercive operator and $b \in V'$. Given a reference control $u_\mathsf d \in H$ and desired state $y_\mathsf d \in H$, we define the multiobjective obstacle problem
\begin{align}
\label{eq:obstacle_mop}
\begin{split}
    &\min\frac{1}{2}\left( \begin{array}{c}
        \lVert y - y_\mathsf d \rVert_H^2  \\[1mm]
        C\lVert u - u_\mathsf d \rVert_H^2
    \end{array} \right), \\
    &\hspace{0,5mm}\text{s.t. }(y, u) \in K \times H\text{ solves } {\langle Ay - u, v - y \rangle}_{V', V} \ge 0 \text{ for all } v \in K,
\end{split}
\end{align}
using a hyperparameter $C > 0$. To reformulate this problem in alignment with the definition of the general MOP \eqref{eq:general_MOP}, we use the control-to-state operator $S:H \to K$ and substitute $y = S(u)$ in the definition of \eqref{eq:obstacle_mop}, which allows to drop the variational inequality.

\begin{figure}
    \centering
    \parbox[b]{0.6\textwidth}{
        \centering 
        \parbox[b]{0.28\textwidth}{
                \centering 
                \includegraphics[width=0.28\textwidth]{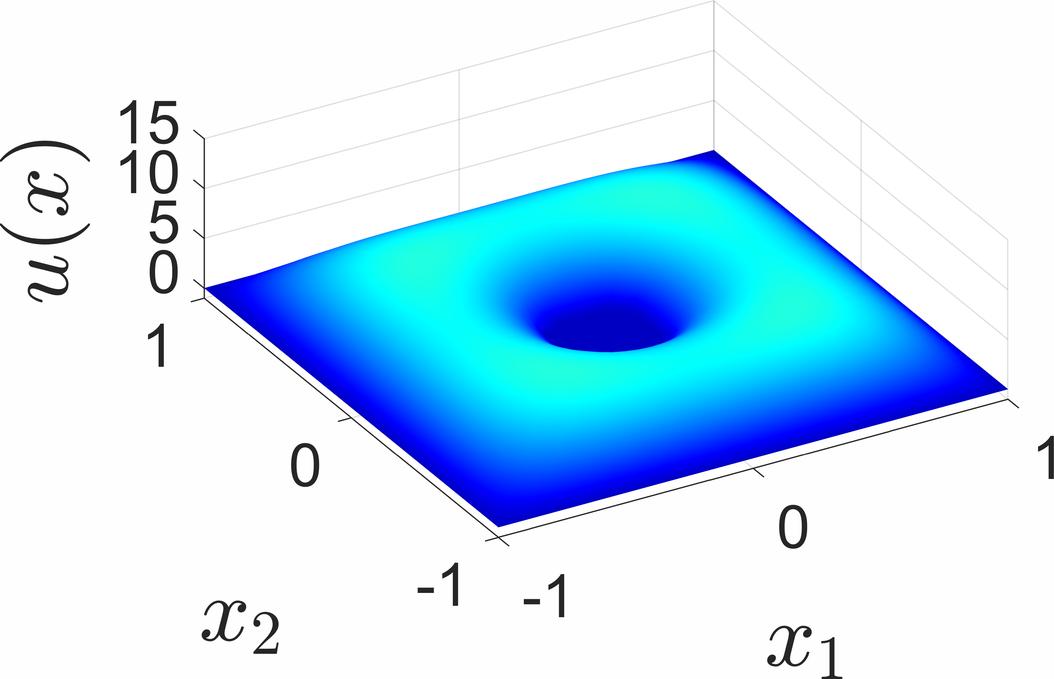}\\
            (a)
    	}
        \parbox[b]{0.28\textwidth}{
            \centering 
            \includegraphics[width=0.28\textwidth]{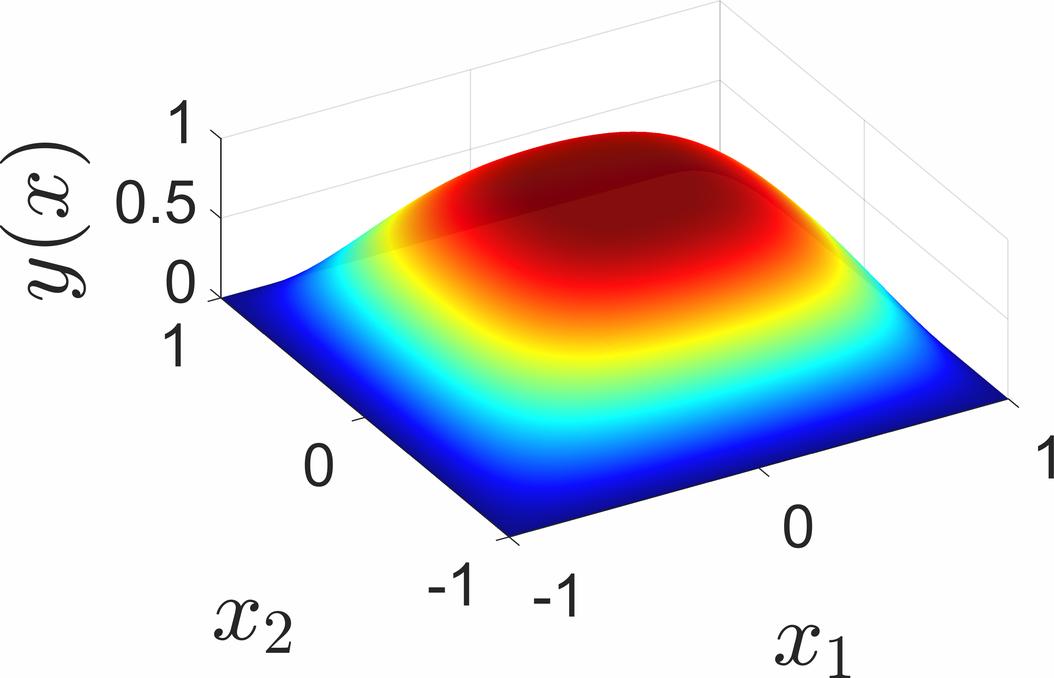}\\
            (b)
            }\\
        \parbox[b]{0.28\textwidth}{
            \centering 
            \includegraphics[width=0.28\textwidth]{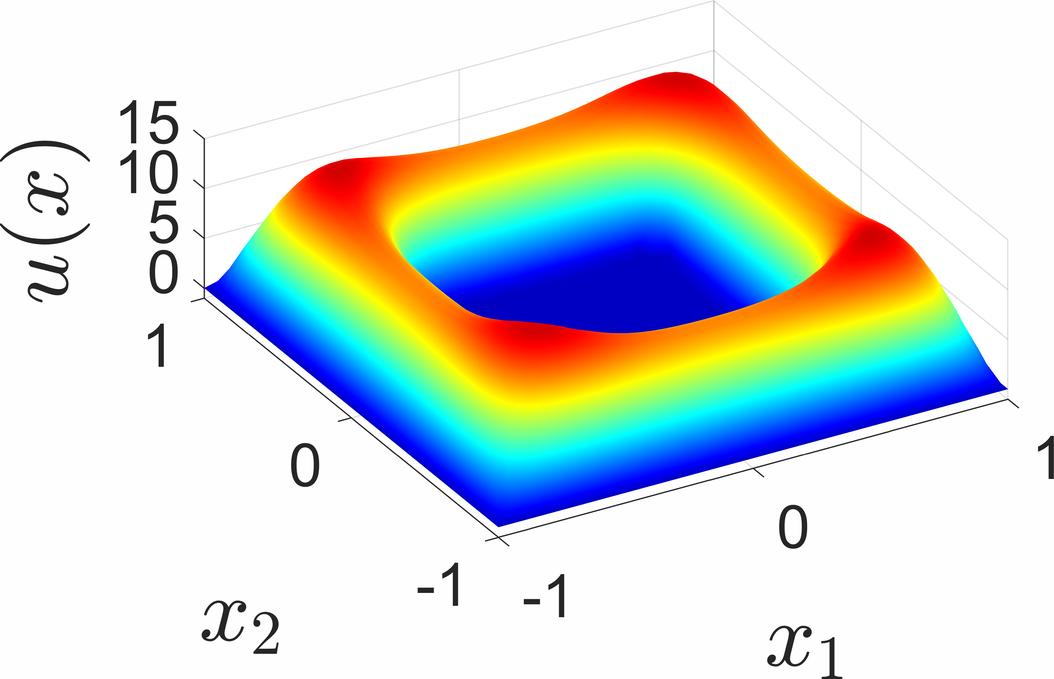}\\
            (c)
    	}
        \parbox[b]{0.28\textwidth}{
            \centering 
            \includegraphics[width=0.28\textwidth]{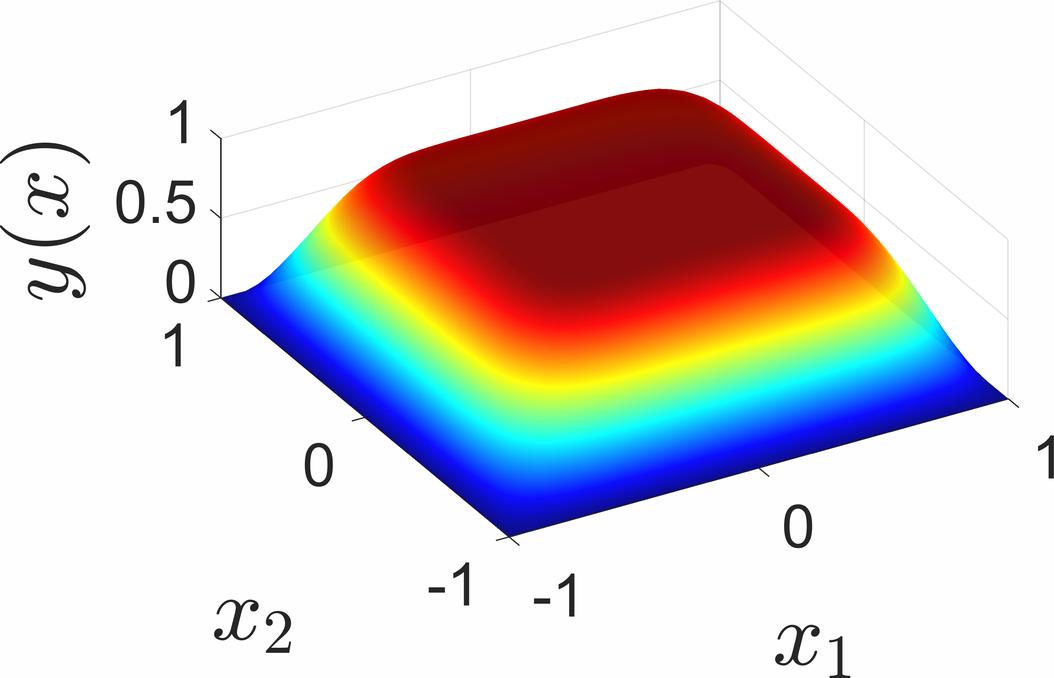}\\
            (d)
    	}
        }
    \parbox[b]{0.39\textwidth}{
        \centering 
        \includegraphics[width=0.39\textwidth]{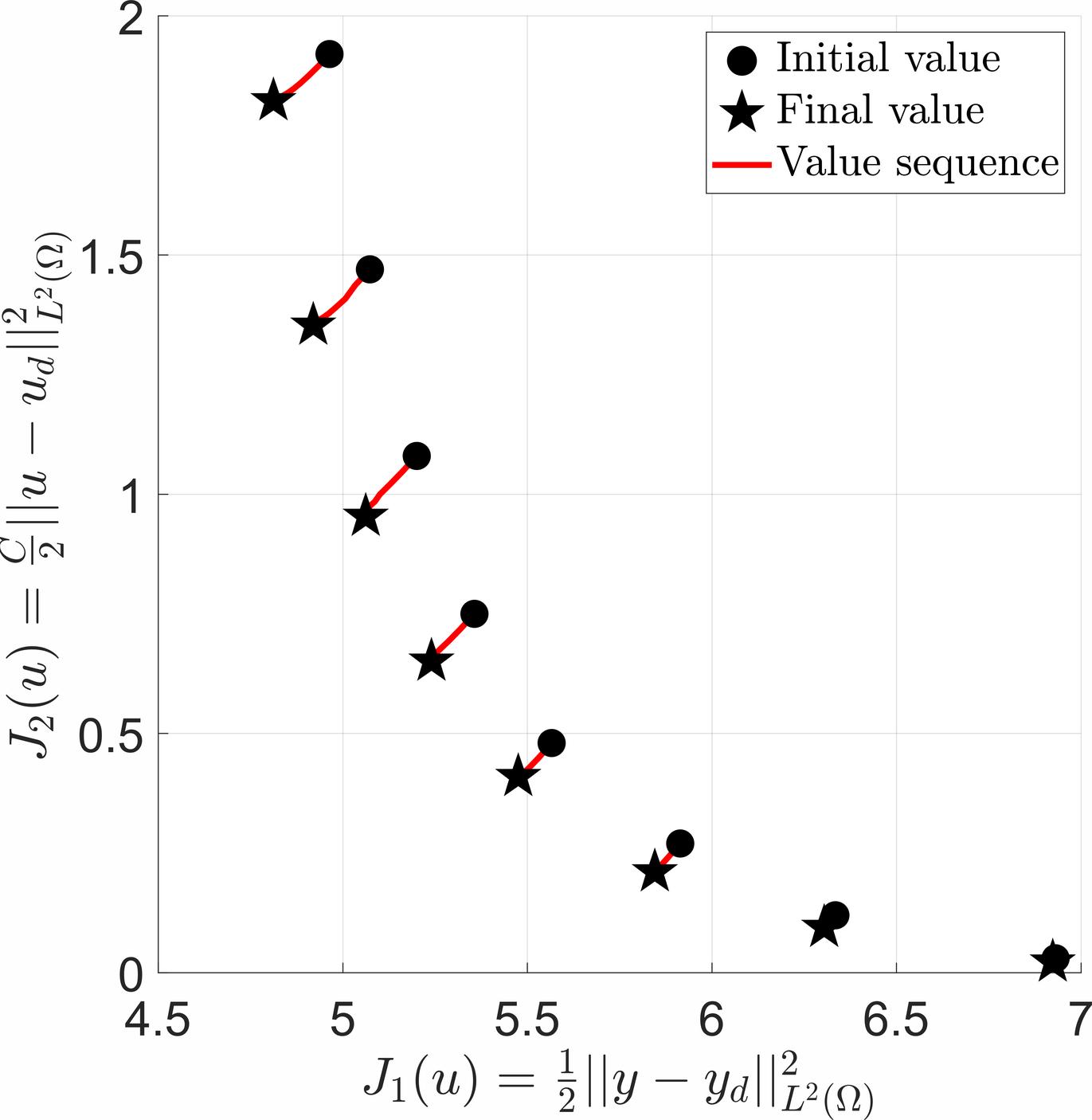}\\
        (e)
        }
    \caption{(a) Pareto optimal control computed from initial control $u \equiv 4$. (b) State corresponding to control from (a). (c) Pareto optimal control computed from initial control $u \equiv 8$. (d) State corresponding to control from (c). (e) Approximation of Pareto front computed from 8 initial values with constant initial controls $u \equiv c$ with $c \in \{1, 2, \dots, 8\}$.}
    \label{fig:example_descent_method_hilbert_space}
\end{figure}

We consider an instance of problem \eqref{eq:obstacle_mop} with reference control $u_d \equiv 0$, desired state $y_d \equiv 2$, constant obstacle $\psi \equiv 1$ and operator $A = -\Delta$, to model a "membrane" $y$ in the plane which gets deformed by a force $u$ and an immovable obstacle $\psi$. By the choice of the reference control and the desired state the objective functions in problem \eqref{eq:obstacle_mop} are conflicting. Using MATLAB's PDE-toolbox we define an FEM-discretization of $\Omega$ to solve problem \eqref{eq:obstacle_mop} numerically. Figure \ref{fig:example_descent_method_hilbert_space} shows parts of the results of an experiment from \cite{SGMPV23}. The figure contains plots of two optimal controls computed with the descent method for different initial values in Subfigures (a) and (c). The corresponding states can be seen in Subfigures (b) and (d). An approximation of the Pareto front computed by starting the multiobjective descent method from multiple initial values is contained in Subfigure (e). For both solutions the state $y$ is bounded from above by the obstacle $\psi \equiv 1$ and we see that the control $u$ vanishes in the areas with contact.


\subsection{Continuation method for sparse optimization}
\label{Section:3.2}

In many applications, one is more interested in finding a simple, good solution of a problem than in finding a complicated, exact solution. For example, when training a neural network, computing the exact minimizer of the loss function leads to overfitting. In many cases, ``simplicity'' of a solution corresponds to sparsity, which can be measured using the $\ell_1$-norm $\| x \|_1 := |x_1| + \dots + |x_n|$ \cite{T1996}. Thus, one is interested in finding solutions with a small $\ell_1$-norm.

Formally, consider a smooth function $L : \R^n \rightarrow \R$. The classical way of computing ``sparse minimizers'' of $L$ is to solve the so-called \emph{regularized problem}
\begin{align}
    \label{eq:sparse_reg_problem}
    \min_{x \in \R^n} L(x) + \lambda\,{\| x \|}_1,
\end{align}
where $\lambda \geq 0$ is the so-called \emph{regularization parameter}. Clearly, the choice of the regularization parameter has a crucial impact on the solution of above problem: If $\lambda$ is chosen too small, then the solution is close to optimal for $L$ but not sparse. If $\lambda$ is chosen too large, than the solution is sparse but does not minimize $L$ at all. Since the ``correct'' $\lambda$ is difficult to choose a priori, one typically computes the solution of \eqref{eq:sparse_reg_problem} for multiple different $\lambda$. In other words, one computes an approximation of the so-called \emph{regularization path}
\begin{align}
    R := \left\{ x \in \R^n : \exists \lambda \geq 0 \text{ with } x \in \argmin_{x \in \R^n} L(x) + \lambda\,{\| x \|}_1 \right\},
\end{align}
which consists of all solutions of \eqref{eq:sparse_reg_problem} for all $\lambda \geq 0$. 

In the context of multiobjective optimization, Problem \eqref{eq:sparse_reg_problem} is a (scaled) weighted sum \cite{Mie12} $\alpha_1 L(x) + \alpha_2 \| x \|_1$ of the non-smooth MOP
\begin{align} \label{eq:reg_MOP}
    \min_{x \in \R^n} 
    \begin{pmatrix}
        L(x) \\
        \| x \|_1
    \end{pmatrix}
\end{align}
for the weights $\alpha_1 = 1/(1+\lambda)$ and $\alpha_2 = \lambda / (1+\lambda)$. In particular, it is easy to see that the regularization path $R$ is contained in the set of weakly Pareto optimal points of \eqref{eq:reg_MOP}. Thus, regularization problems may be solved using results and methods from MOO. Furthermore, since the weighted sum approach is only able to compute the entire Pareto set in case all objectives are convex \cite{Mie12}, treating regularization via MOO may yield more regularized solutions than the classical approach \eqref{eq:sparse_reg_problem}.

In theory, we could use the descent method from Section \ref{Section:3.1} to compute Pareto optimal points of \eqref{eq:reg_MOP}. But since the non-smoothness in \eqref{eq:reg_MOP} occurs in a specific way, it is much more efficient to first exploit the structure of the problem. To this end, in \cite{BGP22}, it was shown that the Pareto critical set $P_c$ of \eqref{eq:reg_MOP} is piecewise-smooth, and that the smooth pieces are Pareto critical sets of MOPs of the form
\begin{align} \label{eq:reg_MOP_proj}
    \min_{x' \in \R^{m}} 
    \begin{pmatrix}
        L(h(x')) \\
        \| x' \|_1
    \end{pmatrix},
    \text{ where }
    (h(x'))_j := 
    \begin{cases}
        x'_i, & j = j_i \\
        0, & \text{otherwise}
    \end{cases}
\end{align}
and $A = \{j_1, \dots, j_m \} \subseteq \{1, \dots, n\}$, $j_1 < \dots < j_m$.
\begin{figure}
    \centering
    \parbox[b]{0.49\textwidth}{
        \centering 
        \includegraphics[width=0.45\textwidth]{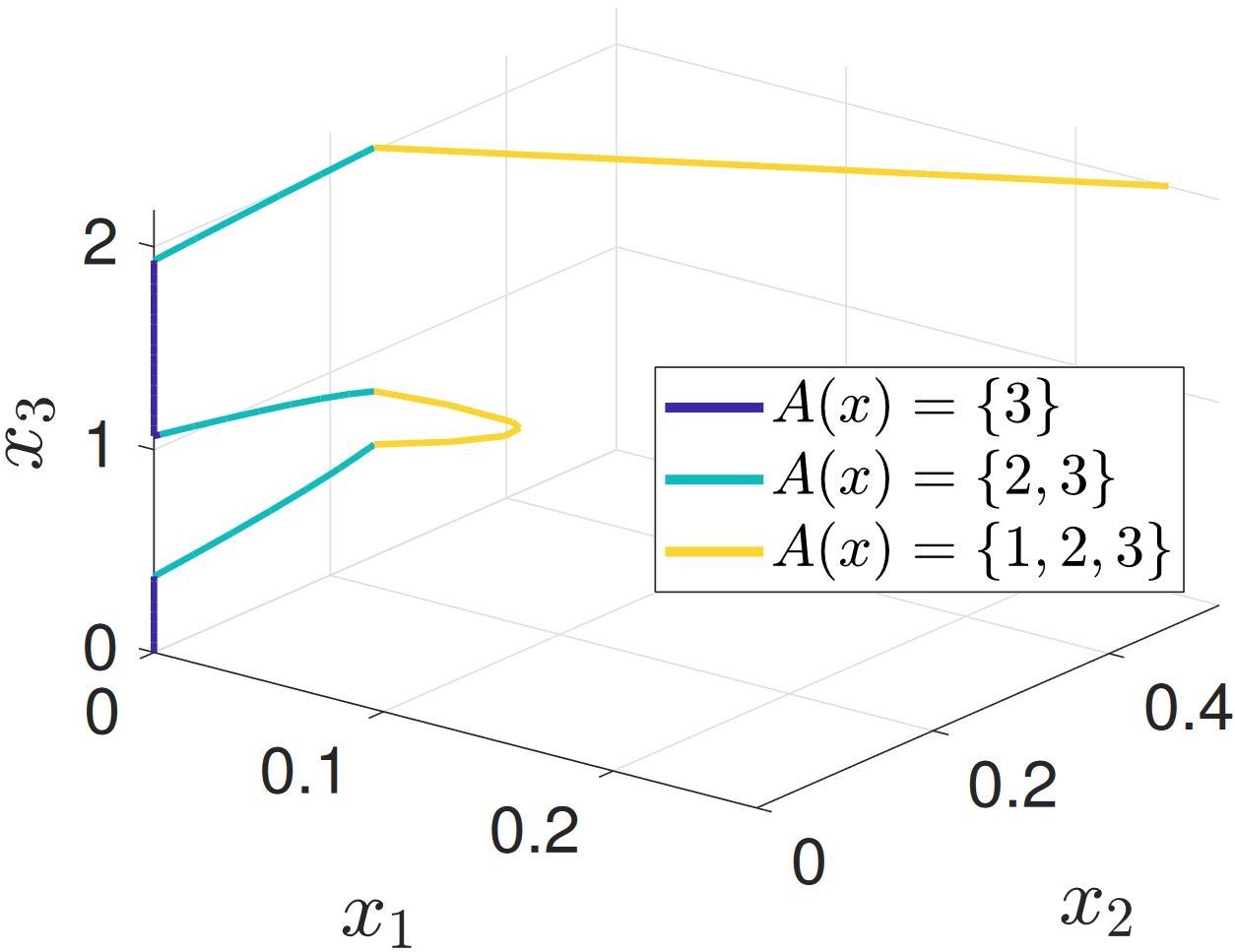}\\
        (a)
	}
    \parbox[b]{0.49\textwidth}{
        \centering 
        \includegraphics[width=0.425\textwidth]{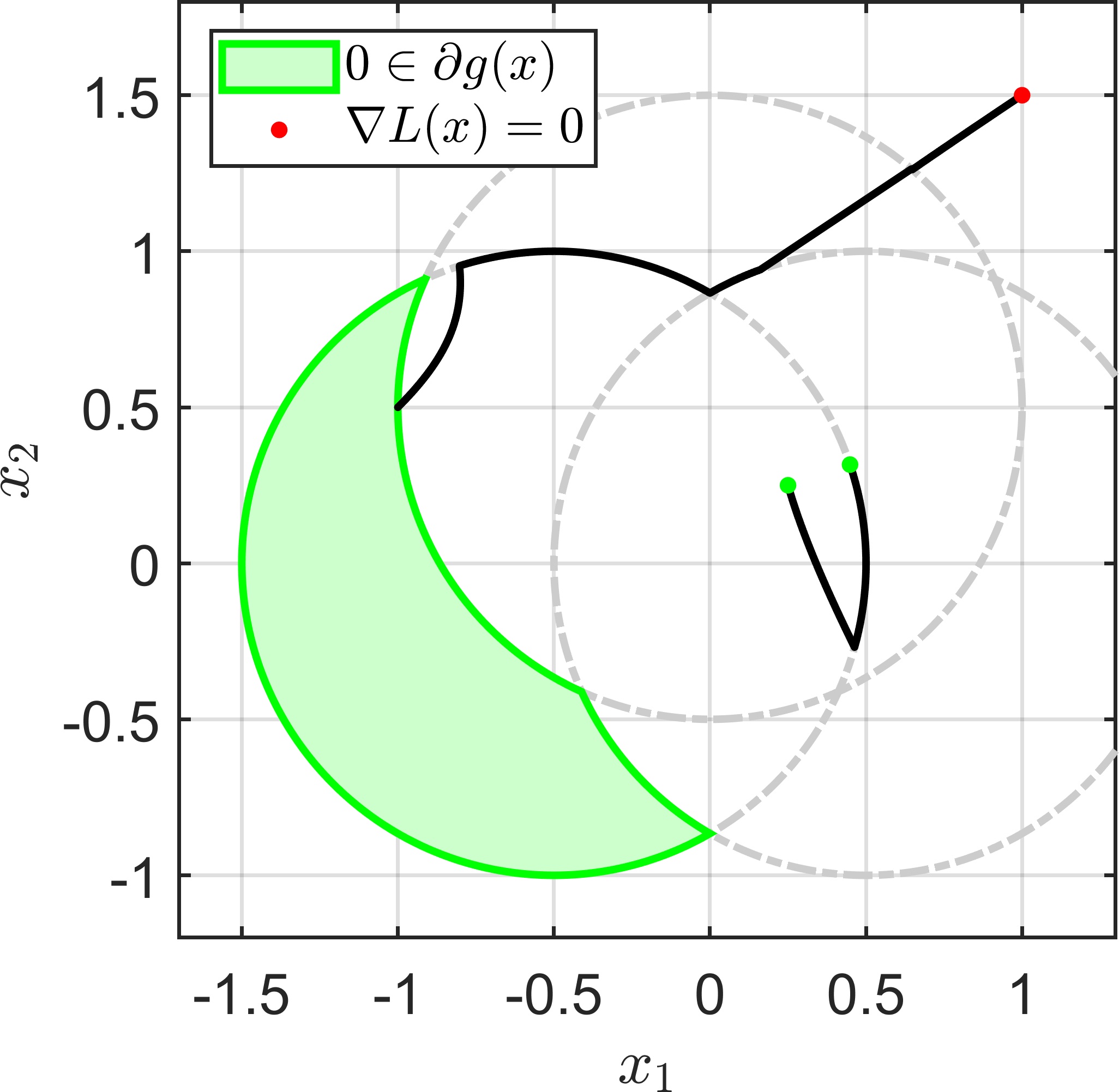}\\
        (b)
	}
    \caption{(a) The decomposition of the Pareto critical set via \eqref{eq:reg_MOP_proj} for the toy example in \cite{BGP22} (Taken from \cite{BGP22}; License: CC BY 4.0). (b) The regularization path (black) for the exact penalty method in Example 6 in \cite{GBP23} (Taken from \cite{GBP23}; License: CC BY 4.0).}
    \label{fig:example_regularization}
\end{figure} 
A visualization is shown in \ref{fig:example_regularization}(a). Since these problems only have to be solved in areas where none of the components of $x'$ vanish, they can be treated and solved as smooth MOPs. As we are aiming at computing a fine approximation of $P_c$, we use \emph{continuation methods} \cite{H2001,SDD2005} for the solution of \eqref{eq:reg_MOP_proj}. The idea of continuation methods is to explore the Pareto critical set from a given starting point, by first computing a \emph{predictor step} along a tangent vector of $P_c$, and then using a local solution method to return to the Pareto critical set in a \emph{corrector step}. This yields a new Pareto critical point close to the initial point, and iteratively repeating this procedure allows us to compute entire connected components. However, as the Pareto critical set of \eqref{eq:reg_MOP} is only piecewise smooth, we need a mechanism to identify kinks in $P_c$ and to figure out the correct set $A$ for problem \eqref{eq:reg_MOP_proj}. In \cite{BGP22}, such a mechanism was derived based on the first-order optimality conditions of \eqref{eq:reg_MOP}, enabling the full exploration of $P_c$ via continuation.

In \cite{GBP23}, the theoretical results from \cite{BGP22} were extended to the case where the second objective in \eqref{eq:reg_MOP} is an arbitrary, twice piecewise differentiable function $g$. This more general setting covers many other regularization problems, e.g., \emph{support-vector machines} \cite{HTF2009}, \emph{total variation denoising} \cite{C2004} and \emph{exact penalty methods} \cite{BKM2014}. In this case, $P_c$ is still piecewise-smooth, but the decomposition into the smooth pieces now also depends on the geometry of the set of non-smooth points of $g$, making it more challenging to construct continuation methods. Figure \ref{fig:example_regularization}(b) shows an example for the regularization path for the exact penalty method from \cite{GBP23}.

\subsection{Other related work}
\label{Section:3.3}

Let us mention some of our related work here:
\begin{itemize}
    \item We derive efficient algorithms to compute weakly Pareto optimal solutions for smooth, convex and unconstrained multiobjective optimization problems in general Hilbert spaces. In \cite{SP22,SP23} a novel inertial gradient-like dynamical system is defied in the multiobjective setting, whose trajectories converge weakly to Pareto optimal solutions.
    \item The hierarchical structure of the Pareto critical sets are studied in \cite{GPD19}.
    \item In \cite{SCM+20} the Pareto Explorer tool is presented—a global/local exploration tool for the treatment of many-objective optimization problems.
\end{itemize}


\section{Surrogate modeling for non-smooth PDE-constrained optimization}
\label{Section:4}

Their are recent results on reduced-order based approaches to MOPs. We mention here \cite{Ban21,BBV16,BBV17,BDPV18,POD19}.

\begin{myexample}
    Let us come back to Example~{\em\ref{Ex:MaxPDE}}. In an RB approach one computes FE solutions to \eqref{MaxPDE3} for suitable parameters $\mu_1,\ldots,\mu_\ell\in\Mad$ with $\ell\ll N$. After a Gram-Schmidt orthogonalization these FE solution span the RB ansatz and test space $V^\ell=\mathrm{span}\,\{\psi_1,\ldots,\psi_\ell\}\subset V^N$ for a Galerkin projection. The choice of the parameters are done by greedy strategies {\em\cite{Haa17,Hesthaven,QMN16}}. Then, we solve the reduced-order problem
\begin{align}
    \label{MaxPDE4}
    \begin{aligned}
         &\int_\Omega c(\mu)\nabla y_\mathsf{rb}(\bx)\cdot\nabla \psi_i(\bx)+a(\mu)\max\big\{0,y_\mathsf{rb}(\bx)\big\}\psi_i(\bx)\,\mathrm d\bx\\
         &\hspace{20mm}=\int_\Omega u(\bx)\psi_i(\bx)\,\mathrm d\bx
    \end{aligned}
\end{align}
for $i=1,\ldots,\ell$ with
\begin{align*}
     y_\mathsf{rb}(\bx)=\sum_{j=1}^\ell\mathrm y_j^\ell\psi_j(\bx)\quad\text{for }\bx\in\Omega
\end{align*} 
with unknowns $\mathrm y_1^\ell,\ldots,\mathrm y_\ell^\ell\in\R$. To be efficient it is necessary to apply the discrete empirical interpolation method (DEIM) to evaluate approximately the non-linear term
\begin{align*}
    \int_\Omega\max\big\{0,y_\mathsf{rb}(\bx)\big\}\psi_i(\bx)\,\mathrm d\bx\quad\text{for }i=1,\ldots,\ell; 
\end{align*}
cf. {\em\cite{CS10}}. Then it turns out that \eqref{MaxPDE4} can be solved much faster than \eqref{MaxPDE3}; cf. {\em\cite{Ber19,BMV22c}}. This allows us to speed-up significantly the numerical optimization. Newly developed a-posteriori error estimators ensure that we rely on surrogate models with sufficient accuracy.
\end{myexample}

To solve \eqref{eq:general_MOP} in the context of PDE constraints a suitable scalarization method in that case is the Pascoletti-Serafini (PS) scalarization \cite{Eic08,PS84}: For a chosen reference point $z\in\R^k$ and a given target direction $r\in\R^k$ with $r > 0$ (component-wise) the Pascoletti-Serafini problem is given by
\begin{align}
    \label{PS-Problem}
    \min \tau\quad\text{s.t.}\quad(\tau,x)\in\R\times\Xad\text{ and }f(u)-z\le \tau r
\end{align}
In \cite{BMV22a} this problem is solved by an augmented Lagrangian approach. However, in our case the evaluation of the objective $f(x)$ requires the solution of a PDE for the given point $x\in\Xad$. This implies further that for the computation of the gradients $\nabla f_i(x)$, $i=1,\ldots,k$ (possibly many adjoint PDEs have to be solved; cf. Example~\ref{Example:PDEs-2}-2). Here, surrogate models offer a promising tool to reduce the computational effort significantly. In an offline phase, a low-dimensional surrogate model of the PDE is constructed by using, e.g., the greedy algorithm, cf. \cite{Haa17,Hesthaven,QMN16}. In the online phase, only the RB model is used to solve the PDE, which saves a lot of computing time. As a new approach \cite{BMV22a} we propose an extension of the method in \cite{Ban21} for solving multi-objective PDE-constained parameter optimization problems. This procedure is based on a combination of a trust-region reduced basis method \cite{BKMOSV22,KMOSV21,QGVW17} and the PS method. In particular, we discuss different strategies to handle the increasing number of reduced basis functions, which is crucial in order to guarantee good performances of the algorithm. We also mention \cite{BP21,BP22} for a trust-region based approach to MOPs.

Let us also refer to the recent work \cite{BV22}, where the following non-smooth equation is considered: For a parameter $\mu \in \Mad$ we consider for $t \in (0,T]$:
\begin{subequations}
    \label{SpaceTime-1}
    \begin{align}
        y_t(t,\bx)-c(\mu)\Delta y(t,\bx)+a(\mu)\max\{0,y(t,\bx)\}&=g(t,\bx)&& \text{for }(t,\bx)\in Q,\\
        y(t,\bx)&=0&&\text{for }(t,\bx)\in\Sigma,\\
        y(0,\bx)&=0&& \text{for }\bx\in\Omega,
    \end{align}
\end{subequations}
where $T$, $\Omega$, $Q$ and $\Sigma$ are as in Example~\ref{Example:PDEs-1}-1). Further, we have $V=H^1_0(\Omega)$ and $H=L^2(\Omega)$. We assume that $c:\Mad\to\R$ is Lipschitz-continuous, positive, uniformly bounded away from zero, that $a:\Mad\to\R$ is Lipschitz-continuous, non-negative and that $f\in C([0,T];H)$ is Lipschitz-continuous. The goal is to derive a space-time formulation of \eqref{SpaceTime-1}. Let us mention that space-time methods have been considered by many authors for (smooth) parabolic problems; see, e.g., to the work \cite{Hin20,LS15,Mei11,Ste15,SY18,Urb14} for (smooth) parabolic problems but there is no error analysis done for nonsmooth PDEs.  To derive the space-time formulation we test \eqref{SpaceTime-1} by test functions $\phi\in \mathscr Y=L^2(0,T;V)$ and integrate over $\Omega$ and $[0,T]$. Then, we derive that the weak solution $y\in\mathscr X=W(0,T)$ to \eqref{SpaceTime-1} satisfies the variational problem
\begin{align}
    \label{SpaceTime-2}
    \mathcal A(y,\phi;\mu) & = \mathcal F(\phi)\quad \text{for all }\phi \in \mathscr Y,
\end{align}
where
\begin{align*}
        \mathcal A(y,\phi;\mu)&=\int_0^T{\langle y_t(t,\cdot),\phi(t,\cdot)\rangle}_{V',V}\,\mathrm dt\\
        &\quad+\int_0^T\int_\Omega\nabla y(t,\bx)\cdot\nabla\phi(t,\bx)+\max\big\{0,y(t,\bx)\big\}\phi(t,\bx)\,\mathrm d\bx\mathrm dt,\\
        \mathcal F(\phi)&=\int_0^T\int_\Omega f(t,\bx)\phi(t,\bx)\,\,\mathrm d\bx\mathrm dt.
\end{align*}
Existence and uniqueness of the solution to (2) follows e.g. from \cite[Theorem 30.A]{Zei89b}.

Analogously to the continuous setting we introduce a discretized space-time formulation. Therefore let $\mathscr X_\delta\subset\mathscr X$ and $\mathscr Y_\delta\subset\mathscr Y$ be finite dimensional subspaces. For $\mu \in \Mad$ we call $y_\delta=y_\delta(\mu)\in \mathscr X_\delta$ a discretized weak solution to \eqref{SpaceTime-1} if
\begin{align}
    \label{SpaceTime-3}
    \mathcal A(y_\delta,\phi_\delta;\mu) & = \mathcal F(\phi_\delta)\quad \text{for all } \phi_\delta \in\mathscr Y_\delta.
\end{align}
In the remainder of this paper we will focus on the case
\begin{align*}
    \mathscr X_\delta=S_{\Delta t}\otimes V_h,\quad \mathscr Y_\delta= Q_{\Delta t}\otimes V_h,
\end{align*}
where $\otimes$ denotes the tensor product and $Q_{\Delta t}$, $S_{\Delta t}$ are piecewise constant, respective piecewise linear finite elements in time and $V_h$ are piecewise linear finite elements in space. As $\delta=(\Delta t, h)$ we summarize the temporal and spatial discretization parameters. 

Analogously the space-time RB setting can be formulated. For a spatial RB space $V_\ell= \text{span}\{\psi_1,\ldots,\psi_\ell\}\subset V_h$ of dimension $\ell\in\mathbb N$, we introduce the RB solution and test spaces
\begin{align*}
    \mathscr X_\mathsf{rb}=S_{\Delta t}\otimes V_\ell,\quad \mathscr Y_\mathsf{rb}= Q_{\Delta t}\otimes V_\ell,
\end{align*}
respectively, where $\mathsf{rb}=(\Delta t,\ell)$ stands for the temporal discretization and RB parameter. For $\mu \in\Mad$ we call $y_\mathsf{rb}=y_\mathsf{rb}(\mu)\in\mathscr X_\mathsf{rb}$ an RB solution to \eqref{SpaceTime-2} if
\begin{align}
    \label{SpaceTime-4}
    \mathcal A(y_\mathsf{rb},\phi_\mathsf{rb};\mu) & = \mathcal F(\phi_\mathsf{rb})\quad \text{for all } \phi_\mathsf{rb} \in\mathscr Y_\mathsf{rb}.
\end{align}

It turns out that \eqref{SpaceTime-3} and \eqref{SpaceTime-4} can be interpreted as a Crank-Nicolson scheme. This can be used for the numerical realization. However, the space-time formulation is used to develop a certified a-posteriori error estimator. This error estimator is adopted to the presence of the DEIM (see \cite{CS10}) as approximation technique for the non-smoothness. The separability of the estimated error into an RB and a DEIM part then guides the development of an adaptive RB-DEIM algorithm, combining both offline phases into one. Numerical experiments show the capabilities of this novel approach in comparison with classical RB and RB-DEIM approaches.

\section{Applications in machine learning}
\label{Section:5}
In this section, we do no longer consider PDE constraints. Instead, we will focus on another application area where highly expensive optimization problems occur very frequently: machine learning.
In supervised machine learning, the central task is basically function approximation. Given input and output tuples, we want to identify a function that maps inputs to output with the smallest possible error. Put differently, we need to choose the model's parameters in such a way that we minimize the \emph{empirical loss}, see \cite{AML12,GBC16} for a detailed introduction. 

If the number of parameters becomes large (as in deep learning) or the data set grows in size, these training problems tend to become very expensive to solve.
Traditionally, one thus only considers a single objective during training. However, we can easily identify multiple criteria such as the empirical loss, model sparsity, the inclusion of physical knowledge \cite{KKL+21}, or a measure of interpretability, to name just a few. It is thus of vital importance to improve the performance of multicriteria methods in order to make them usable in the context of machine learning, which is currently a highly active area of research.

In the following, we will address two exemplary applications where data science and multiobjective optimization meet. In Section \ref{Section:5.1}, this concerns the identification of conflicting criteria from a given data set.

\subsection{Inferring the objectives of an MOP from data}
\label{Section:5.1}
In \emph{inverse optimization}, the goal is to infer the objective function (or parameters that influence it) from its optimizers. In scalar optimization, the solution of a problem is generically a single point, not restricting the possible objective functions much. In multiobjective optimization on the other hand, the Pareto set typically consists of infinitely many solutions, containing much more information than just a single point and allowing us to infer more about the objectives.

Formally, assume that we are given a finite set of points $D \subseteq \R^n \times \R^k$ consisting of pairs $(x,\alpha) \in D$. Our goal is to find a continuously differentiable ($C^1$) function $f : \R^n \rightarrow \R^k$ such that for all $(x,\alpha) \in D$, $x$ is Pareto critical for $f$ with multiplier $\alpha$ as in \eqref{eq:smooth_KKT}. To this end, we choose a set of basis functions $b_1, \dots, b_d \in C^1(\R^n,\R)$. Then for $f_i = \sum_{j = 1}^d c_{ij} b_j$ and fixed $x$ and $\alpha$, the convex combination in \eqref{eq:smooth_KKT} is linear in $c \in \R^{k \cdot d}$. Thus, considering the optimality conditions for all $(x,d) \in D$ yields a linear system $\mathcal{L} c = 0$. If this system has a nontrivial solution $c^*$, then setting $f_i = \sum_{j = 1}^d c^*_{ij} b_j$ yields a desired objective vector.
More generally, in \cite{GP21b}, it was shown that if $c^*$ is a right-singular vector corresponding to the smallest singular value $s$ of $\mathcal{L}$, then
\begin{align*}
    \left\| \sum_{i = 1}^k \alpha_i \nabla f_i(x) \right\|_2 \leq s \quad \forall (x,\alpha) \in D
\end{align*}
for such $f$. Furthermore, numerical experiments were carried out, showing potential applications for test problem generation, stochastic MOO and surrogate modelling.


\section{Mixed-integer techniques for efficient data-driven surrogate modeling of PDE-constrained control problems}
\label{Section:6}
In this final section, we are treating a quite different class of problems. Instead of having non-smoothness in our problem statement, we here treat continuous control problems, but we manually introduce non-smoothness for efficiency purposes.
As this is a new concept, we begin with the single objective setting. However, we the extension to multiple objectives can be realized using different approaches for multiobjective optimal control, see, e.g., \cite{PD18a}.

As motivated above, the central task here is to efficiently solve \emph{optimal control problems} for smooth but complex and thus expensive-to-evaluate dynamical systems.
Mathematically speaking, we consider the following problem over the time horizon $p\cdot \Delta t$:
\begin{equation} \tag{I} \label{eq:OCP}
\begin{aligned}
\min_{u\in U^p} J(y) &= \min_{u\in U^p} \sum_{i=0}^{p-1} P(y_{i+1}) \\
\mbox{s.t.} \quad y_{i+1} &= \Phi(y_i, u_i), \qquad i = 0,1,2,\ldots,
\end{aligned}
\end{equation}
where $y_i$ and $u_i \in U$ are the system state and control at time instant $t_i = i \Delta t$, with $U$ being the set of admissible controls, e.g., $U=[u^{\mathsf{min}}, u^{\mathsf{max}}]$. 
The objective function (for instance, the distance to some desired trajectory $y^{\mathsf{ref}}$) is denoted by $P$, and $\Phi$ describes the flow of the underlying dynamical system (e.g., an ordinary or a partial differential equation) over the time increment $\Delta t$. 
The solution of \eqref{eq:OCP} yields the optimal control $u^*$ and corresponding state $y^*$. Feedback control can then be achieved via \emph{Model Predictive Control (MPC)} \cite{GP17}, i.e., by solving \eqref{eq:OCP} repeatedly over a short horizon and applying the first entry $u^*_0$ to the real system which is running simultaneously. Here we mentioned related work on MPC approaches for MOPs in \cite{BFGV20,HPW+18,HPW+19,HOP20,OP21}.

A substantial challenge that we often face is the fact that the efficient prediction (and, by extension, control) of complex dynamical systems is hindered by the fact that the system dynamics are either very expensive to simulate or even unknown. Researchers have been investigating ways to accelerate the solution by using data for decades, the \emph{Proper Orthogonal Decomposition (POD)} being an early and very prominent example \cite{Sir87}. 
More recently, the major advances in data science and machine learning have lead to a plethora of new possibilities, for instance artificial neural networks such as \emph{Long Short-Term Memory (LSTM)} Networks \cite{HS97} or \emph{Reservoir Computers / Echo State Networks} \cite{JH04}, regression-based frameworks for the identification of nonlinear dynamics \cite{BPK16}, or numerical approximations of the \emph{Koopman operator} \cite{KNP+20,PK19,RMB+09,Sch10}, which describes the linear dynamics of observable functions. These methods facilitate the efficient simulation and prediction of high-dimensional spatio-temporal dynamics using measurement data, without requiring prior system knowledge. 
As a consequence of the success of data-driven prediction, many approaches for control have been presented over the past decades. 
However, a drawback is that the construction of surrogate models with inputs is often much more tedious and also problem-specific and data hungry \cite{BCB05,BPB+20}.

\begin{figure*}[b]
	\centering
	\includegraphics[width=.8\textwidth]{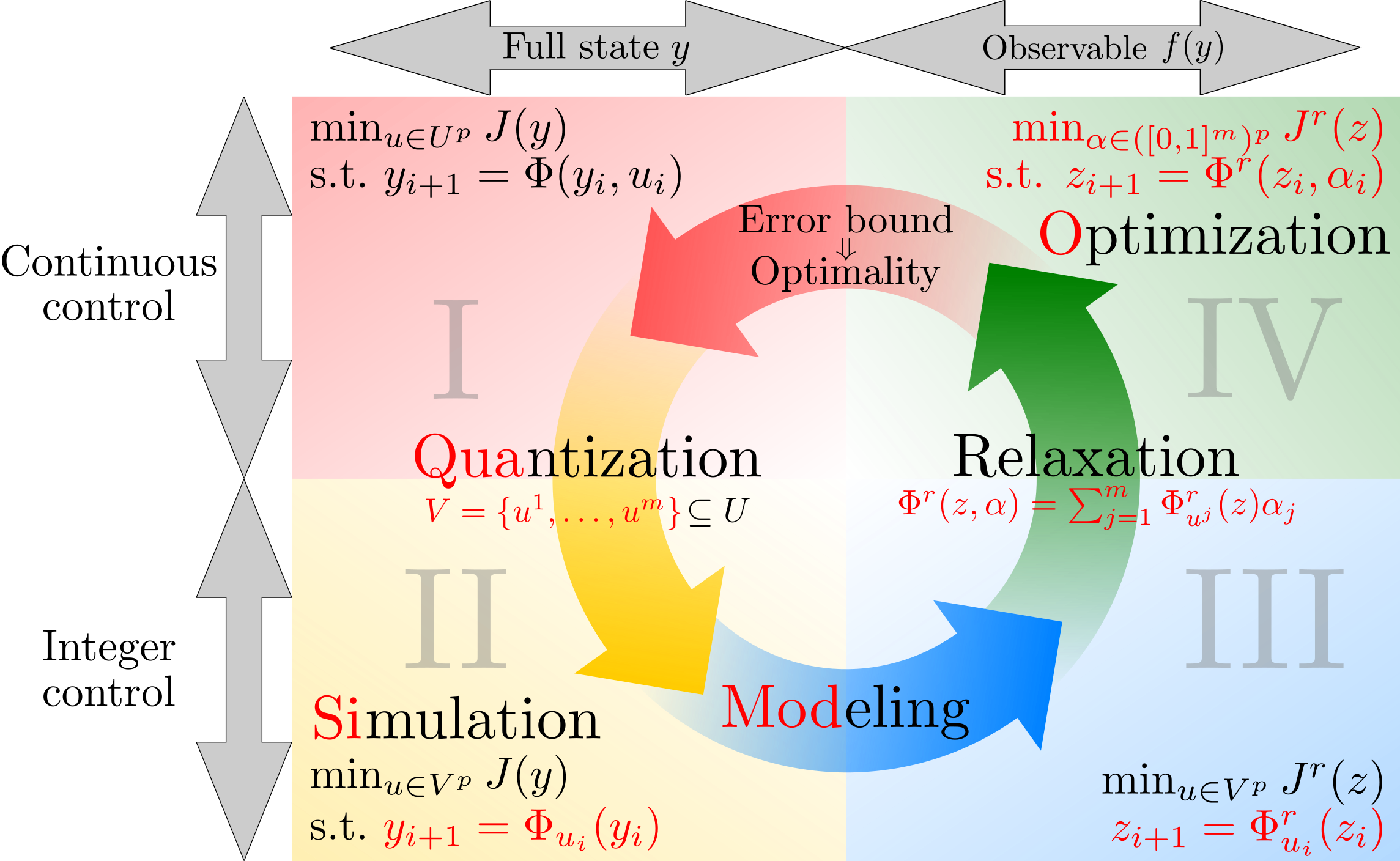}
	\caption{The QuaSiModO framework consisting of the four steps Quantization, Simulation, Modeling and Optimization (Taken from \cite{PB23}; License: CC BY 4.0).}
	\label{fig:QuaSiModO}
\end{figure*}

The approach we present here to solve \eqref{eq:OCP} via surrogate models while avoiding the aforementioned issues is based on modifying the control problem instead of adjusting the surrogate modeling to the control setting.
The resulting framework, which we call \emph{QuaSiModO}, consists of the following steps (cf.~also Figure \ref{fig:QuaSiModO}):
\begin{enumerate}
	\item {\textbf{Qua}ntization} of the the admissible control $U$ (for instance by replacing the interval $U=[u^{\mathsf{min}}, u^{\mathsf{max}}]$ by the bounds $V=\{u^{\mathsf{min}}, u^{\mathsf{max}}\}$);
	\item {\textbf{Si}mulation} of the autonomous systems (e.g., $\Phi_{u^{\mathsf{min/max}}}(y) = \Phi(y,{u^{\mathsf{min/max}}})$;
	\item {\textbf{Mod}eling} of the individual systems -- using either the full state $y$ or some observable $z=f(y)$ -- via an arbitrary ``off-the-shelf'' surrogate modeling technique (POD, neural network, Koopman operator, etc.);
	\item {\textbf{O}ptimization} using the resulting set of autonomous surrogate models and relaxation techniques.
\end{enumerate}
This interplay between continuous and integer control modeling as well as between the full system state and observed quantities (e.g., measurements) allows us to utilize the best of both worlds, namely
\begin{itemize}
	\item integer controls for efficient data-driven modeling using arbitrary predictive models,
	\item continuous control inputs for real-time control, and
	\item existing error bounds for predictive models.
\end{itemize}
QuaSiModO successively transforms Problem \eqref{eq:OCP} into related control problems that -- as long as the predictive surrogate model is sufficiently accurate -- yield optimal trajectories $y^*$ that are close to one another. From (I) to (II), we quantize the control, meaning that only a finite set $V \subseteq U$ of inputs is feasible. This allows us to replace the non-autonomous dynamical system $\Phi(y, u)$ by a finite set of autonomous systems $\Phi_{u^j}(y)$, each corresponding to one entry $u^j \in V$. While introducing an artificial drawback from the control perspective (Problem (II) is a mixed-integer optimal control problem), we can now easily introduce an equivalent Problem (III) that is based on surrogate models $\Phi^r_{u^i}(z)$ for a reduced quantity $z = f(y)$. Here, the function $f$ is an \emph{observable} which maps measurements from the state space of the full system to the space of measurements (which may be of significantly smaller dimension). 
As the transformation from (II) to (III) acts on a set of autonomous systems, we can approximate the individual systems $\Phi_{u^j}$ from individual measurement data sets, using whichever method we prefer. 

In order to mitigate the disadvantages with respect to the complexity of the control problem, the problem of selecting an optimal input from $V$ is relaxed by determining the optimal convex combination of the autonomous systems:
\begin{equation} \tag{IV} \label{eq:rOCP}
    \begin{aligned}
        \min_{\alpha\in ([0,1]^m)^p} J^r(z) &= \min_{\alpha\in ([0,1]^m)^p} \sum_{i=0}^{p-1} P^{r}(z_{i+1}) \\
        \mbox{s.t.} \qquad~~ z_{i+1} = \Phi^r(z_i, \alpha_i) &= \sum_{j=1}^m \alpha_{i,j} \Phi^r_{u^j}(z_i) \quad \mbox{and} \quad
        \sum_{j=1}^m \alpha_{i,j} = 1.
    \end{aligned}
\end{equation}
Problem \eqref{eq:rOCP} is again continuous -- with respect to the input $\alpha$. 
For control affine systems, we can now determine $u^* = \sum_{j=1}^m \alpha^*_j u^j$ and directly apply it to the real system. For non-affine systems, we use the sum up rounding algorithm from \cite{SBD12}, by which a control corresponding to one of the quantized inputs is applied to the real system. 

\begin{figure}[b!]
	\centering
	\includegraphics[width=\textwidth]{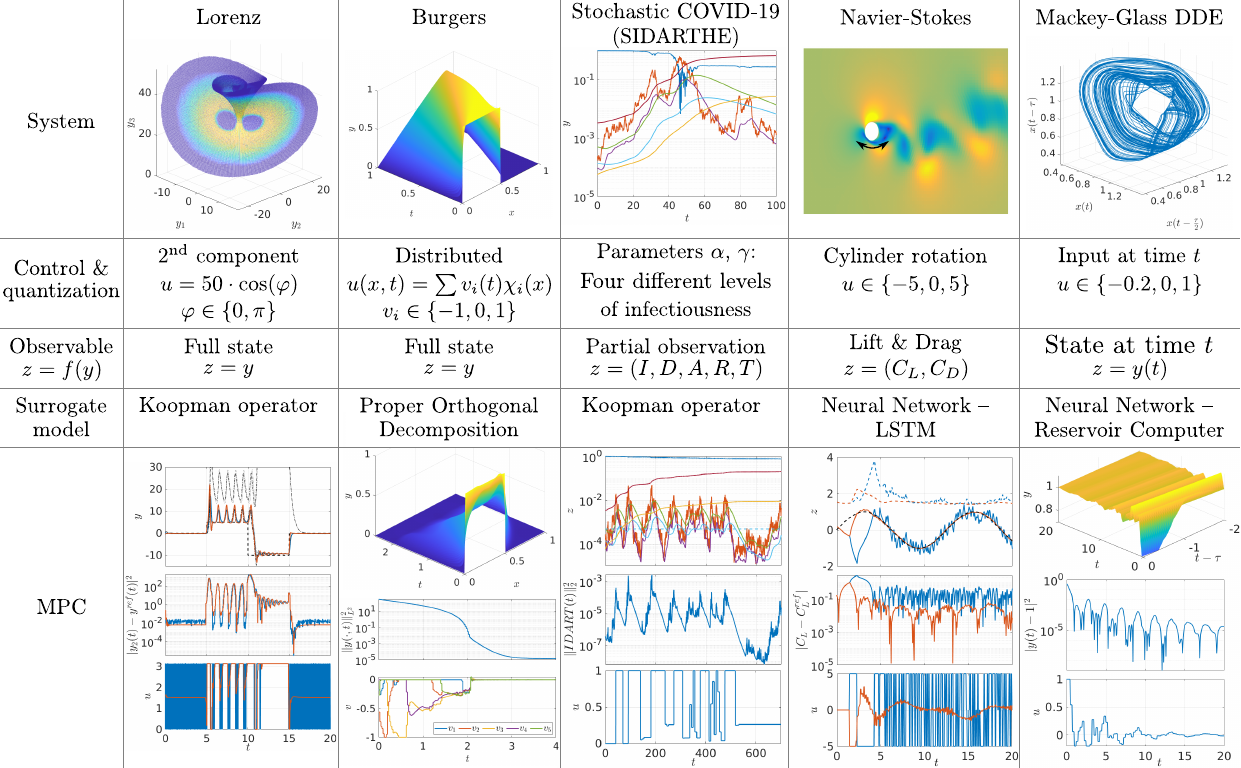}
	\caption{MPC using QuaSiModO applied to various combinations of systems and surrogate models (Taken from \cite{PB23}; License: CC BY 4.0).}
	\label{fig:Results}
\end{figure}

Besides the ability to include arbitrary predictive models into the QuaSiModO framework, an important aspect is that existing error bounds for the chosen surrogate model can easily be included, see \cite{PB23} for a detailed description
The availability of error bounds is of particular importance for engineering systems, where safety is of utmost importance (e.g., for aircraft or autonomous vehicles). The bounds guarantee the performance of a controller and -- more importantly -- will automatically become stronger with future developments in the field of data-driven modeling.

We have tested the QuaSiModO framework on a variety of dynamical systems, observable functions and surrogate modeling techniques, cf.~Figure~\ref{fig:Results}, a detailed description is given in \cite{PB23}. 
For instance, we can control the lift force acting on a cylinder (determined by the velocity and pressure fields governed by the 2D Navier--Stokes equations) without any knowledge of the flow field using the standard LSTM framework included in \emph{TensorFlow},
and stabilize the Mackey-Glass equation using a standard echo state network.
This highlights the flexibility and broad applicability of the method and the success of the technique in constructing data-driven feedback controllers.



\bibliographystyle{plain}
\bibliography{literature}

\end{document}